\documentclass[12pt]{amsart}
\usepackage{latexsym}           
\usepackage{amssymb}
\usepackage{epsfig}             
\usepackage{amsfonts}

\newcommand\beq{\begin{equation}}
\newcommand\eeq{\end{equation}}

\newcommand{\IP}{\mathbb{P}}                                     
\newcommand{\IR}{\mathbb{R}}                           
\newcommand{\IC}{\mathbb{C}}
\newcommand{\IZ}{\mathbb{Z}}

\newcommand{\M}{\mathcal{M}}

\newcommand{\cP}{\mathcal{P}}

\newcommand{\cB}{\mathcal{B}}
\newcommand{\cC}{\mathcal{C}}

\newcommand{\cF}{\mathcal{F}}

\newcommand{\cO}{O}            

\newcommand{\h}{\mathfrak{h}}

\newcommand{\lsl}{\mathfrak{sl}} 

\newcommand{\bt}{{\bf t}}
\newcommand{\bm}{{\bf m}}

\newcommand{\pf}{\begin{bpf}}

\newcommand{\pfms}{\begin{bpfms}}
\newcommand{\epf}{\end{bpf}\hfill$\square$\\}           
\newcommand{\epfms}{\end{bpfms}\hfill$\square$\\}               

\newcommand{\wt}{\widetilde}

\newcommand{\wh}{\widehat}
\newcommand{\al}{\alpha}
\newcommand{\be}{\beta}
\newcommand{\ga}{\gamma}
\newcommand{\de}{\delta}

\newcommand{\si}{\sigma}

\newcommand{\ze}{\zeta}

\newcommand{\bL}{{\bf\Lambda}}

\newcommand{\ad}{\text{\rm ad}}

\newcommand{\re}{{\text{\rm Re}}}

\newcommand{\tr}{\text{\rm Tr}}
\newcommand{\Hom}{\text{\rm Hom}}

\newcommand{\SL}{\text{\rm SL}}
\newcommand{\PSL}{\text{\rm PSL}}
\newcommand{\GL}{\text{\rm GL}}
\newcommand{\PGL}{\text{\rm PGL}}
\newcommand{\SU}{\text{\rm SU}}

\newcommand{\End}{\text{\rm End}}
\newcommand{\diag}{{\text{\rm diag}}}

\newcommand {\epi}{\varepsilon_\infty}
\newcommand {\epii}{\varepsilon^{-1}_\infty}
\newcommand {\eps}{\varepsilon}
\newcommand {\epsi}{\varepsilon^{-1}}

\newcommand{\spq}{/\!\!/}

\def\mapright#1{\smash{
        \mathop{\longrightarrow}\limits^{#1}}}

\def\hookmapright#1{\smash{
        \mathop{\hookrightarrow}\limits^{#1}}}

\theoremstyle{plain}
\newtheorem {hypo}{\bf\hspace{-\parindent}Hypothesis}

\newtheorem {thm}{Theorem}
\newtheorem {prop}[hypo]{Proposition}

\newtheorem {cor}[hypo]{Corollary}

\newtheorem {lem}[hypo]{Lemma}

\theoremstyle{definition}

\theoremstyle{remark}
\newtheorem {rmk}[hypo]{Remark}

\begin{document}



\title[From Klein to Painlev\'e]
{From Klein to Painlev\'e \\ via Fourier, Laplace and Jimbo}
\author{Philip Boalch}
\address{\'Ecole Normale Sup\'erieure\\
45 rue d'Ulm\\
75005 Paris\\
France} 
\email{boalch@dma.ens.fr}

\begin{abstract}
We will describe a method for constructing explicit algebraic solutions to the
sixth Painlev\'e equation, generalising that of Dubrovin--Mazzocco.
There are basically two steps: 
First we explain how to construct finite braid group orbits of triples
of elements of $\SL_2(\IC)$  out of triples of generators of three-dimensional 
complex reflection groups.
(This involves the Fourier--Laplace transform for certain irregular  
connections.) 
Then we adapt a result of Jimbo to produce the Painlev\'e VI solutions.
(In particular this solves a Riemann--Hilbert problem explicitly.)

Each step will be illustrated using the complex reflection group associated to
Klein's simple group of order 168.
This leads to a new algebraic solution with seven branches.
We will also prove that, unlike the algebraic solutions of Dubrovin--Mazzocco
and Hitchin, this solution is not equivalent to any
solution coming from a finite subgroup of $\SL_2(\IC)$.

The results of this paper also yield a simple proof of a
recent theorem of Inaba--Iwasaki--Saito 
on the action of Okamoto's affine $D_4$ symmetry group as well as the correct
connection formulae for generic Painlev\'e VI equations.
\end{abstract}


\maketitle


\renewcommand{\baselinestretch}{1.02}            
\normalsize

\begin{section}{Introduction}

Klein's quartic curve
$$X^3Y+Y^3Z+Z^3X = 0\quad\subset\ \IP^2(\IC)$$
is of genus three and 
has the maximum possible number $84(g-1)=168$ 
of holomorphic automorphisms.
Klein found these automorphisms explicitly
(in terms of $3\times 3$ matrices). They constitute {Klein's simple group}
$K \subset \PGL_3(\IC)$ 
which is isomorphic to $\PSL_2(7)$.

Lifting to $\GL_3(\IC)$ 
there is a two-fold covering group $\wh K\subset \GL_3(\IC)$
of order $336$ which is a complex reflection group---there 
are complex reflections
\beq\label{klein gens}
r_1,\ r_2,\ r_3\ \in \GL_3(\IC)
\eeq
which generate $\wh K$.
(Recall a pseudo-reflection is an automorphism 
of the form ``one plus rank one'', 
a complex reflection is a pseudo-reflection of finite order and 
a complex reflection group is a finite group generated by 
complex reflections. Here, each generator $r_i$ 
has order two---as for real reflections.) 

Using the general tools to be described in this paper we will construct, 
starting from the Klein complex reflection group $\wh K$,
another algebraic curve with affine equation
\beq\label{eq: soln curve intro}
F(t,y)=0
\eeq
given by a polynomial $F$ with integer coefficients.
This curve will be a seven-fold cover of the $t$-line
branched only at $0,1,\infty$ and such that the function
$y(t)$, defined implicitly by \eqref{eq: soln curve intro}, 
solves the
Painlev\'e VI differential equation.

One upshot of this will be to construct an explicit 
rank three Fuchsian 
system of linear differential equations with four singularities
(at $0,t,1,\infty$, for some $t$) on $\IP^1$, and with monodromy group
equal to $\wh K$ in its natural representation (so the monodromy around each of
the 
finite singularities $0,t,1$ is a generating reflection).

In general the construction of linear differential equations with finite
monodromy group 
is reasonably straightforward provided one works
with rigid representations of the monodromy groups.
In our situation the representation is minimally non-rigid; it lives in a
complex two-dimensional moduli space, and this is the basic reason 
the (second order) Painlev\'e VI equation arises.

Apart from the many physical applications,
from a mathematical perspective our main interest in the Painlev\'e 
VI equation is that it is the explicit form of the simplest 
isomonodromy (=non-abelian Gauss--Manin) connection.
In brief, the isomonodromy connections arise by replacing
the closed differential forms and periods appearing
in the usual (abelian) Gauss--Manin picture, by flat connections and 
monodromy
representations, respectively. 

Indeed one may view the Painlev\'e VI equation as a natural nonlinear
analogue of the Gauss hypergeometric equation. 
From this point of view the thrust of this paper is towards 
finding the analogue of Schwartz's famous list of hypergeometric equations
with algebraic solutions.

Before carefully describing the contents of this paper we will 
briefly recall exactly how the sixth Painlev\'e equation arises.

Consider a Fuchsian system of
differential equations (with four singularities) 
of the form
\begin{equation} \label{eqn: 2x2 linear}
\frac{d\Phi}{dz}= A(z)\Phi; \qquad A(z)=\sum_{i=1}^3\frac{A_i}{z-a_i}
\end{equation}
where the $A_i$'s are $2\times 2$ traceless matrices. 
We wish to deform \eqref{eqn: 2x2 linear} isomonodromically---i.e.
when the 
pole positions $(a_1,a_2,a_3)$ are moved in 
$\IC^3\setminus{\text{diagonals}}$ we wish to 
vary the coefficients $A_i$ such that
the conjugacy class of the corresponding monodromy representation
is preserved.
Such isomonodromic deformations are governed by Schlesinger's equations:
\begin{equation} \label{eqn: 2x2 schles}
\frac{\partial A_i}{\partial a_j}=\frac{[A_i,A_j]}{a_i-a_j}\qquad 
\text{if } i\ne j, \text{ and}\qquad
\frac{\partial A_i}{\partial a_i}=-\sum_{j\ne
i}\frac{[A_i,A_j]}{a_i-a_j}.
\end{equation}

Let us view these more geometrically as a nonlinear connection 
on a fibre-bundle.
First observe that Schlesinger's equations preserve the adjoint orbit $O_i$
containing each $A_i$ and 
are invariant under overall conjugation of $(A_1,A_2,A_3,A_4)$,
where $A_4=-A_1-A_2-A_3$ is the residue of \eqref{eqn: 2x2 linear} at infinity.
Thus one sees 
that Schlesinger's equations amount to a flat connection,
the {\em isomonodromy connection}, on the
trivial fibre bundle 
\begin{equation}\label{eqn: rk 2+ bundle}
\M^*:=(O_1\times O_2\times O_3\times O_4)\spq G \times B \mapright{} B
\end{equation}
over $B:=\IC^3\setminus{\text{diagonals}}$, where the fibre 
$(O_1\times\cdots \times O_4)\spq G$ is the quotient of 
$$
\left\{\ (A_1,A_2,A_3,A_4)\in O_1\times O_2\times O_3\times O_4 
\ \bigl\vert \ \text{$\sum A_i =0$}\ \right\}$$ 
by overall conjugation by $G=\SL_2(\IC)$. 
(Generically this fibre is two dimensional 
and has a natural complex symplectic structure.) 

Now for each point $(a_1,a_2,a_3)$ of the base $B$ one can also 
consider the set 
\beq\label{eq: monod fibres}
\Hom_\cC(\pi_1(\IC\setminus\{a_i\}),G)/G
\eeq
of conjugacy classes of 
representations of the 
fundamental group of the four-punctured sphere, where the representations are
restricted to take the simple loop around $a_i$
into the conjugacy class $\cC_i:=\exp(2\pi\sqrt{-1} \cO_i)\subset G$ 
($i=1,\ldots,4, a_4=\infty$).
These spaces  of representations are also generically two dimensional (and
complex symplectic) and fit
together into a fibre bundle
$$M\mapright{} B.$$
Moreover this bundle $M$ has a {\em complete} flat connection  
defined locally by identifying representations
taking the same values on a fixed  set of fundamental group generators.
The isomonodromy connection is the pullback of this complete
connection 
along the natural bundle map
$$\nu: \M^*\mapright{}M$$ 
defined by 
taking the systems \eqref{eqn: 2x2 linear} to their monodromy representations
(cf. \cite{Hit95long, smid}).

To obtain Painlev\'e VI one 
chooses specific coordinates $x, y$ on the fibres of 
$\M^*$ and then, 
upon restricting the pole positions to $(a_1,a_2,a_3)=(0,t,1)$,
the isomonodromy connection amounts to 
two first order coupled nonlinear equations for $x(t), y(t)$.
Eliminating $x$ yields (cf. e.g. \cite{JM81}) 
the sixth Painlev\'e equation\footnote
{The general PVI equation was first written down by R. Fuchs (son of
L. Fuchs) and it was added to the list of Painlev\'e equations 
by Painlev\'e's student B. Gambier.}  
(PVI):

$$\frac{d^2y}{dt^2}=
\frac{1}{2}\left(\frac{1}{y}+\frac{1}{y-1}+\frac{1}{y-t}\right)
\left(\frac{dy}{dt}\right)^2
-\left(\frac{1}{t}+\frac{1}{t-1}+\frac{1}{y-t}\right)\frac{dy}{dt}  $$
$$
\quad\ +\frac{y(y-1)(y-t)}{t^2(t-1)^2}\left(
\al+\be\frac{t}{y^2} + \gamma\frac{(t-1)}{(y-1)^2}+
\delta\frac{ t(t-1)}{(y-t)^2}\right). $$

The four parameters $\al,\be,\ga,\delta\in\IC$  here are directly 
related to the choice of the adjoint orbits $O_i$. 
From another viewpoint, we will see the monodromy spaces
\eqref{eq: monod fibres} are affine cubic surfaces 
and Iwasaki
\cite{Iwas-modular}
has recently
pointed out that the four parameters
correspond to the moduli of such surfaces, appearing in the Cayley normal
form.

The sixth Painlev\'e equation
has critical singularities at $0,1,\infty$ and 
is remarkable in that any of its solutions have 
wonderful analytic continuation properties: any locally-defined 
solution $y(t)$  may be analytically continued to a meromorphic function on
the universal cover of the three-punctured sphere
$\IP^1\setminus\{0,1,\infty\}$. (This is the so-called Painlev\'e property.)

From the geometric viewpoint, the monodromy of PVI (i.e. the analytic
continuation of solutions around $\IP^1\setminus\{0,1,\infty\}$)
corresponds to the monodromy of the nonlinear connection on $\M^*$.
In turn this connection is the pullback of the complete connection on the
bundle $M$.
Being complete, the monodromy of the connection on $M$ amounts to an action
of the fundamental group of the base $B$ (the pure braid group $\cP_3$)
on the standard 
fibre \eqref{eq: monod fibres}.
This is the standard braid group action on the monodromy data,
which we thus see gives the monodromy of solutions 
to PVI. \footnote{
Restricting to 
$\IP^1\setminus\{0,1,\infty\}\hookmapright{\iota}B$,
where $\iota(t)=(0,t,1)$, 
amounts to restricting to the action of the free subgroup 
$\cF_2:= \pi_1(\IP^1\setminus\{0,1,\infty\})\hookmapright{\iota_*} \cP_3$
of the braid group.
This $\cF_2$ action is equivalent to the $\cP_3$ action since
the centre $Z\cong \IZ$ of $\cP_3$ acts trivially 
and $\iota_*$ induces an isomorphism $\cF_2\cong \cP_3/Z$.}

Our main concern in this paper is to construct algebraic solutions 
to PVI.
One knows that, for generic values of the four parameters, any solution of
PVI is a `new transcendental function' on the universal cover of the
three-punctured sphere.
However, for special values of the parameters it is possible that there 
are solutions expressible in terms of standard transcendental functions,
or even solutions which are algebraic---i.e. are defined by polynomial
equations.
For example there are the algebraic 
solutions of Hitchin \cite{Hit-Poncelet, Hit-Octa}, Dubrovin \cite{Dub95long}
and  Dubrovin--Mazzocco \cite{DubMaz00} 
related to the dihedral,
tetrahedral, octahedral and icosahedral groups.

The problem of constructing  algebraic solutions may be broken into two
steps.
First, the algebraic solutions will have a finite number of branches
and so one  may start by looking for finite orbits of the $\cP_3$ action on the
space of monodromy data \eqref{eq: monod fibres}.

Clearly if we start with a linear system \eqref{eqn: 2x2 linear}
whose monodromy is a finite subgroup of $\SL_2(\IC)$, then the corresponding 
braid group orbit will be finite.
The solutions of Hitchin, Dubrovin and  
Mazzocco mentioned above are equivalent to
solutions arising in this way.

The basic idea underlying the present paper is that PVI also arises 
as the equation for isomonodromic deformations of certain rank three 
Fuchsian systems. 
Namely we replace $A_1,A_2,A_3$ in \eqref{eqn: 2x2 linear}
by $3\times 3$ matrices $B_1,B_2,B_3$ each of rank one.
Then the corresponding moduli spaces are still of dimension two,
and one finds again that PVI governs the isomonodromic deformations
(and that any PVI equation arises in this way).
Note that the rank one condition implies the 
corresponding monodromy group will be generated by a triple of
pseudo-reflections in $\GL_3(\IC)$.

Rather than work throughout with this equivalent $3\times 3$ representation
of PVI, we will pass between the two pictures in order to use 
existing machinery developed in the $2\times 2$ framework
(in particular the work of Jimbo \cite{Jimbo82}).

Our starting point will be to describe a method of constructing 
finite braid group orbits of triples of elements of $\SL_2(\IC)$
starting from any triple of complex reflections generating 
a complex reflection group in $\GL_3(\IC)$.
In general this will yield more exotic finite braid group orbits
than those from finite subgroups of $\SL_2(\IC)$.
The key idea behind this construction is to use the Fourier--Laplace
transformation to convert the rank three Fuchsian system into a 
rank three system with an irregular singularity, 
then to apply a simple scalar shift and transform back, 
so that the resulting Fuchsian system is reducible,
and we take the irreducible rank two quotient or subsystem.
Of crucial importance here is Balser--Jurkat--Lutz's computation 
\cite{BJL81} of the action of the Fourier--Laplace
transformation on monodromy data, 
relating the monodromy data of the Fuchsian system to the 
Stokes data $u_\pm$ of the irregular system.
This correspondence may be described by the 
explicit formula
$$r_3r_2r_1=u_-^{-1}t^2u_+$$
(dating back at least to Killing) for the Birkhoff 
factorisation of the product of generating reflections,
and enables us to compute the action of the scalar shift on the reflections.

Lots of finite braid group orbits of $\SL_2(\IC)$ triples are obtained in this
way: 
We recall that Shephard--Todd \cite{Shep-Todd} 
have classified all the complex reflection groups and showed
that in three-dimensions, apart from the real reflection groups, there are four
irreducible complex reflection groups generated by triples of reflections,
of orders 336, 648, 1296 and 2160
respectively, as well as two infinite families $G(m,p,3), m\ge 3, p=1,m$ of
groups of orders $6m^3/p$. 
For $m=2$ and $p=1,2$ 
these would be the symmetry groups
of the octahedron and tetrahedron respectively.
(In general, 
for other $p$ dividing $m$, $G(m,p,3)$ is not generated by a triple of
reflections.)
The main example we will focus on, 
the Klein group, is thus the smallest non-real
exceptional complex reflection group. This leads to a $\cP_3$ orbit of size
seven which
we will prove is not isomorphic to any orbit
coming from a finite subgroup of $\SL_2(\IC)$.

The second step in the construction of algebraic solutions is to pass from the
finite braid group orbit to the explicit solution.
For this we adapt (and correct) a 
result of Jimbo \cite{Jimbo82} giving an explicit
formula for the leading term in the asymptotic 
expansion at zero of the solution $y(t)$
on each branch.
By using the PVI equation this is sufficient to determine the solution curve
precisely:
\begin{gather}
F(t,y)  = \quad \notag \\
\left( 162\,{t}^{3}-243\,{t}^{2}-243\,t+162 \right) {y}^{7}+
\left( -567\,{t}^{3}+2268\,{t}^{2}-567\,t \right) {y}^{6}+ \notag \\
\left( -1701\,{t}^{3}-1701\,{t}^{2} \right) {y}^{5}+ 
\left( 1407\,{t}^{4}+2856\,{t}^{3}+1407\,{t}^{2} \right) {y}^{4}+
\label{soln curve} \\
\left( 14\,{t}^{5}-2849\,{t}^{4}-2849\,{t}^{3}+14\,{t}^{2} \right) {y}^{3}+
\left( -21\,{t}^{5}+3444\,{t}^{4}-21\,{t}^{3} \right) {y}^{2}+\notag \\
\left( -567\,{t}^{5}-567\,{t}^{4} \right) y+
(125\,{t}^{6}-88\,{t}^{5}+125\,{t}^{4})\notag
\end{gather}

\noindent
which admits the rational parameterisation:
\begin{equation}\notag
y=-{\frac { \left( 5\,{s}^{2}-8\,s+5 \right)  \left( 7\,{s}^{2}-7\,s+4\right) }
{ s \left( s-2 \right)  \left( s+1 \right) \left( 2\,s-1 \right) 
\left( 4\,{s}^{2}-7\,s+7 \right)  }},\qquad
t={\frac { \left( 7\,{s}^{2}-7\,s+4 \right) ^{2}}{{s}^{3} \left( 4\,{s}^
{2}-7\,s+7 \right) ^{2}}}.
\end{equation}

Using this parameterisation it is easy to carry out the ultimate test
and substitute back into the Painlev\'e VI equation, with 
$(\al,\be,\ga,\de)=(9,-4,4,45)/98$,
finding that we do indeed have a solution.

Note that we 
have not considered the further problem of writing down the affine Weyl
group orbit of \eqref{soln curve}, or of finding the simplest representative.

The general strategy of this paper is the same as the paper 
\cite{DubMaz00} of Dubrovin--Mazzocco.
Indeed part of our motivation was
to extend their work to (a dense open subset of) the full four parameter
family of PVI equations.
Recall that \cite{DubMaz00} dealt with the real (orthogonal) 
three-dimensional reflection
groups and for this it was sufficient to only consider
a one-parameter family of PVI
equations (corresponding to fixing each of $A_1,A_2,A_3$ to be nilpotent,
so the remaining parameter is the choice of orbit of $A_4$).

In relation to \cite{DubMaz00} the key results of the present paper are
firstly to see how to extend their method of passing from 
generating triples of orthogonal reflections to finite $\cP_3$ 
orbits of (unipotent) $\SL_2(\IC)$ triples. 
(Reading the earlier papers 
\cite{Dub95long,DubPT} of Dubrovin was helpful to 
fully understand this aspect of \cite{DubMaz00}.)

Secondly we were able to fix Jimbo's asymptotic formula.
(Dubrovin--Mazzocco did not use Jimbo's asymptotic result, but adapted 
Jimbo's argument to prove
a version of it for their nilpotent situation.)
The key point here was to find a sign error hidden in the depths of
Jimbo's asymptotic formula---perhaps we should emphasize that 
without the correction the construction of this paper
will not work at all. (Namely at some point 
we need to obtain precise rational numbers out of the transcendental formulae.)
This sign is also important because
it is needed to obtain the correct connection formulae for solutions 
of the Painlev\'e VI equation---by symmetry analogous asymptotic formulae may
be obtained at one and at infinity, thereby giving the connection formulae.

The two main tools of this paper (construction of finite $\cP_3$
orbits
of $\SL_2(\IC)$ triples, and Jimbo's formula)
are independent and will have separate applications. 
For example one may take any triple of elements of a finite subgroup of 
$\SL_2(\IC)$ and try to apply Jimbo's formula to find solutions to PVI.
(E.g. in \cite{icos} we have classified the inequivalent $\cP_3$ 
orbits of generators of the binary icosahedral group 
and, as a further test of Jimbo's formula, constructed 
a new algebraic solution to PVI with 12 branches---this 
is the largest genus zero icosahedral solution and is
interesting since its parameters lie on {\em none} 
of the reflecting hyperplanes
of Okamoto's affine $F_4$ action.)  

\vspace{.4cm}

The layout of this paper is as follows.
In section \ref{sn: braid orbits} we explain in a direct algebraic fashion
how to obtain finite 
braid group orbits of (conjugacy classes of) triples of elements of
$\SL_2(\IC)$ from triples of generators of three-dimensional
complex reflection groups.
Section \ref{sn: imds} (which could be skipped on a first reading)
then explains how the formulae of section \ref{sn: braid orbits} were found.
This is somewhat more technical, involving the action of the 
Fourier--Laplace transform on monodromy data, 
but is necessary to understand 
the origin of the procedure of section \ref{sn: braid orbits}.
We also mention in passing (Remark \ref{rmk: qwgp})  
the relation with the
$\GL_n(\IC)$ quantum Weyl group actions.  
Next, in section \ref{sn: jimbo}, we give Jimbo's formula 
for the leading term in the asymptotic expansion at zero 
of the PVI solution $y(t)$ on the branch specified by a given $\SL_2(\IC)$
triple. 
This is applied in section \ref{sn: klein soln} to find the 
Klein solution explicitly.  
Section \ref{sn: ineq} then proves that the Klein solution is not 
equivalent (under Okamoto's affine $F_4$ action) 
to any solution coming from a finite subgroup of $\SL_2(\IC)$.
Then in section \ref{sn: reconstruction} we explain
how to reconstruct, from such a PVI solution, an explicit rank three Fuchsian 
system with monodromy group generated
by the triple of complex reflections we started with in 
section \ref{sn: braid orbits}.
Finally in section \ref{sn: 3x3 rep} we describe a direct path from the
$3\times 3$ Fuchsian isomonodromic deformations to PVI and deduce a recent 
theorem of Inaba--Iwasaki--Saito \cite{IIS}.

It should be mentioned that, in the short paper \cite{pecr}, 
we previously showed by a different method that the 
equations for 
isomonodromic deformations of the $3\times 3$ Fuchsian systems mentioned above
are equivalent to PVI---this 
was written before Jimbo's formula was fixed and also does not give the
relation between the rank two and three monodromy data.

\ 

\renewcommand{\baselinestretch}{0.95}           
\small

{\small
{\bf Acknowledgments.}\ \ 
The Klein solution was found whilst the author 
was a J. F. Ritt assistant professor
of Mathematics 
at Columbia University, New York, and was announced
at the April 2003 AMS meeting 
at the Courant Institute.
Various other 
parts were done whilst the author was a member of 
DPMMS (Cambridge), The Mathematical Institute (Oxford), SISSA (Trieste) and 
IRMA (Strasbourg, supported by the EDGE Research Training Network
HPRN-CT-2000-00101).
The author is grateful to Nigel Hitchin, Boris Dubrovin and Marta Mazzocco
for useful conversations, and to Anton Alekseev for inviting him to visit
the University of Geneva Mathematics  Department and  the 
Erwin Schr\"odinger Institute (Vienna) during
the summer of 2003, where this paper
was written (supported by the Swiss NSF and the ESI
respectively).
}

\end{section}

\begin{section}{Braid group orbits}\label{sn: braid orbits}

In this section we will explain how to obtain some interesting finite braid
group orbits of triples of elements of $\SL_2(\IC)$ from triples of 
generators of three-dimensional complex reflection groups. 

The motivation is simply the fact that branches of a solution to PVI are
parameterised by pure braid group orbits of conjugacy classes of triples of
elements of $\SL_2(\IC)$.
Clearly any algebraic solution of PVI has a finite number of branches and so
the first step towards finding an algebraic solution is to find a finite braid
group orbit, which is what we will do here.

\subsection*{$\bf2\times 2$ case} \ 

Let $G=\SL_2(\IC)$ and consider the standard 
action of the three-string braid group
$\cB_3$ on $G^3$ generated by
\beq\label{eq: sl2 b3 action}
\be_1(M_3,M_2,M_1)=(M_2,M_2^{-1}M_3M_2,M_1)
\eeq
$$\be_2(M_3,M_2,M_1)=(M_3,M_1,M_1^{-1}M_2M_1)$$
where $M_i\in G$.
We are interested in the induced action of $\cB_3$ on the set of conjugacy
classes of such triples.

First we recall some basic facts (cf. e.g. \cite{Magnus}). 
To begin with note that
the seven functions
$$
m_1:=\tr(M_1),\qquad
m_2:=\tr(M_2),\qquad
m_3:=\tr(M_3),$$
\beq\label{eq: 7 2x2fns}
m_{12}:=\tr(M_1M_2),\qquad
m_{23}:=\tr(M_2M_3),\qquad
m_{13}:=\tr(M_1M_3)
\eeq
$$m_{321}:=\tr(M_3M_2M_1).$$
on $G^3$ are invariant under the diagonal conjugation action of $G$ and in
fact generate the ring of invariant polynomials.
Indeed the ring of invariants is isomorphic to the quotient 
of $\IC[m_1,m_2,m_3,m_{12},m_{23},m_{13},m_{321}]$ by the (ideal generated by
the) so-called Fricke relation:
\beq\label{eq: fricke}
m_{321}^2-Pm_{321}+Q=4
\eeq
where $P,Q$ are the following polynomials in the first six variables:
$$P=m_1m_{23}+m_2m_{13}+m_3m_{12}-m_1m_2m_3$$
$$Q=m_1^2+m_2^2+m_3^2+m_{12}^2+m_{23}^2+m_{13}^2+m_{12}m_{23}m_{13}
-m_1m_2m_{12}-m_2m_3m_{23}-m_1m_3m_{13}.$$
(This appears in the book \cite{FrickeKleinI} of Fricke and Klein.)
That there is precisely one relation fits 
nicely with the rough dimension count of $3\cdot 3-3=6$ for the 
space of conjugacy classes of triples.
Viewed as a quadratic equation for $m_{321}$ the other root of 
\eqref{eq: fricke} is $\tr(M_1M_2M_3)$ so in particular we have
\beq\label{eq 123 2x2}
\tr(M_1M_2M_3)=P-m_{321}.
\eeq
Note that,  upon fixing $m_1,m_2,m_3,m_{321}$, the Fricke relation is a cubic
equation in the remaining three variables;
the six dimensional variety we are studying is essentially a universal family
of affine cubic surfaces \cite{Iwas-modular}.

Now, the induced $\cB_3$ 
action on conjugacy classes of triples induces an action on the invariant
functions, and we will describe this action in terms of 
the seven chosen generators.
Clearly $m_{321}$ is fixed by both $\be_i$, and the $m_i$ are just permuted:
$$
\be_1(m_1,m_2,m_3)=(m_1,m_3,m_2),\qquad
\be_2(m_1,m_2,m_3)=(m_2,m_1,m_3).$$
\begin{lem}
The induced $\cB_3$ action on the quadratic functions is
$$\be_1(m_{12},m_{23},m_{13})=
(m_2m_{321}+m_1m_3-m_{13}-m_{12}m_{23},m_{23},m_{12})$$
$$\be_2(m_{12},m_{23},m_{13})=
(m_{12},m_{13},m_1m_{321}+m_2m_3-m_{23}-m_{13}m_{12})$$
\end{lem}
\pf
The second formula follows from the first by permuting indices.
For the first formula the hard part is to establish
$$\tr(M_2^{-1}M_3M_2M_1)=m_2m_{321}+m_1m_3-m_{13}-m_{12}m_{23}.$$
One way to do this (which will extend to the $3\times 3$ case below)
is to write $M_i=\varepsilon_i(1+e_i\otimes\al_i)$ for some 
rank one matrix $e_i\otimes\al_i$ and number $\varepsilon_i\in\IC^*$.
Then $\tr(M_2^{-1}M_3M_2M_1)$ 
can be expanded in terms of the numbers $\al_i(e_j)$ 
and the terms of the resulting expression can be identified as terms in the
expansions of the seven invariant functions.
\epf

Before moving on to the higher rank case we point out
the evident fact that
if $(M_3,M_2,M_1)$ are a triple of generators of a finite subgroup of 
$G$ then the corresponding braid group orbit is finite 
(and in turn the induced action on conjugacy classes of triples 
is also finite). 

\subsection*{$\bf3\times 3$ case} \ 

Now we wish to find analogous formulae for the corresponding action of $\cB_3$
on conjugacy classes of triples of pseudo-reflections in $\GL_3(\IC)$.

Suppose $r_1,r_2,r_3$ are pseudo-reflections in $\GL_3(\IC)$, so that 
$$r_i= 1+e_i\otimes\al_i$$
for some $e_i\in V,\al_i\in V^*$ where $V=\IC^3$.
Choose six non-zero complex numbers $n_1,n_2,n_3,t_1,t_2,t_3$ such
that $t_i$ is a choice of square root of $\det(r_i)$
(i.e. $t_i^2=1+\al_i(e_i)$), 
that the product $r_3r_2r_1$ has eigenvalues $\{n_1^2,n_2^2,n_3^2\}$
and that these square roots are chosen so that
\beq\label{eq: sqrt of det}
t_1t_2t_3=n_1n_2n_3
\eeq
(which is a square root of the equation $\prod(\det r_i)=\det r_3r_2r_1$).
These square roots (and  the choice of ordering of eigenvalues of $r_3r_2r_1$) 
will not be needed to describe the braid group actions here,
but will be convenient later.

Now consider the following eight $\GL_3(\IC)$-invariant 
functions on the set of triples
of pseudo-reflections:
$$t_1^2,\qquad t_2^2,\qquad t_3^2,$$ 
\beq\label{eqn: t defn}
t_{12}:=\tr(r_1r_2)-1,\qquad 
t_{23}:=\tr(r_2r_3)-1,\qquad 
t_{13}:=\tr(r_1r_3)-1,
\eeq
$$t_{321}:=n_1^2+n_2^2+n_3^2,\qquad
t_{321}':=(n_1n_2)^2+(n_2n_3)^2+(n_1n_3)^2.
$$
(Note that $t_i^2=\tr(r_i)-2,\  
t_{321}:=\tr(r_3r_2r_1)$ and $t_{321}'=\det(r_3r_2r_1)\tr(r_3r_2r_1)^{-1}$.)
The subtractions of $1$ or $2$ in this definition 
turn out to simplify the formulae below.
The action of $\cB_3$ on triples of pseudo-reflections is generated by
\beq\label{eq: b3 on prs}
\be_1(r_3,r_2,r_1)=(r_2,r_2^{-1}r_3r_2,r_1),
\eeq
$$\be_2(r_3,r_2,r_1)=(r_3,r_1,r_1^{-1}r_2r_1).$$

Now consider the induced action on conjugacy classes of triples.
First, it is clear that $t_{321},t_{321}'$ are 
$\cB_3$-invariant since $r_3r_2r_1$ is fixed. 
Also, as before, the functions $t^2_i$ are just permuted:
$$
\be_1(t^2_1,t^2_2,t^2_3)=(t^2_1,t^2_3,t^2_2),\qquad
\be_2(t^2_1,t^2_2,t^2_3)=(t^2_2,t^2_1,t^2_3).$$
\begin{lem} \label{lem: 3x3 quad action}
The induced $\cB_3$ action on the functions 
$(t_{12},t_{23},t_{13})$ is as follows:
$$\be_1(t_{12},t_{23},t_{13})=
(
t_{321}+t^2_1+t^2_3-t_{13}
+(t_{321}'-t_{12}t_{23})/t^2_2
,t_{23},t_{12})$$
$$\be_2(t_{12},t_{23},t_{13})=
(t_{13},t_{23},
t_{321}+t^2_2+t^2_3-t_{23}
+(t_{321}'-t_{13}t_{12})/t^2_1).
$$
\end{lem}
\pf
The non-obvious part is to establish
$$\tr(r_2^{-1}r_3r_2r_1)=t_{321}+t^2_1+t^2_3-t_{13}
+(t_{321}'-t_{12}t_{23})/t^2_2.$$
For this we first observe $r_i^{-1}=1-e_i\otimes\al_i/t_i^2$.
Then expanding  $\tr(r_2^{-1}r_3r_2r_1)$ yields
$$t_{321}+1-t_2^2-u_{12}u_{21}-u_{23}u_{32}
-(u_{12}u_{23}u_{31}+u_{23}u_{32}u_{12}u_{21})/t^2_2$$
where $u_{ij}:=\al_i(e_j)$.
To simplify 
this we first use the following identities (obtained by expanding the
traces $t_{ij}$):
\beq\label{eq uu and t's}
u_{ij}u_{ji}=t_{ij}-t_i^2-t_j^2\qquad\text{if $i\ne j$}.
\eeq
Then, to finish, we use the identity (analogous to \eqref{eq 123 2x2}):
$$u_{12}u_{23}u_{31}=
t^2_3t_{12}+
t^2_2t_{13}+
t^2_1t_{23}
-(t_1t_2)^2-(t_2t_3)^2-(t_1t_3)^2 -t_{321}',$$
which is obtained by expanding 
$\tr(r_1^{-1}r_2^{-1}r_3^{-1})=
\tr\left((r_3r_2r_1)^{-1}\right)=n_1^{-2}+n_2^{-2}+n_3^{-2}.$
\epf

Again we have the evident fact that
if $(r_3,r_2,r_1)$ are a triple of generators of a finite subgroup of 
$\GL_3(\IC)$, i.e. if they are generators of a three-dimensional
complex reflection group,
then the corresponding braid group orbit is finite 
(and in turn the induced action on conjugacy classes of triples 
is also finite). 

\begin{rmk}
A rough dimension count 
gives $3\cdot 5-8=7$ for the space of conjugacy classes 
of pseudo-reflections, so we expect there to be a relation amongst the
eight invariant functions. This is
the analogue of the Fricke relation and comes from the
identity   
$$(u_{12}u_{23}u_{31})(u_{32}u_{21}u_{13})=
(u_{23}u_{32})(u_{12}u_{21})(u_{13}u_{31}).$$
Rewriting each bracketed term in terms of the eight functions
yields the desired relation:
\begin{align} \label{eq: 3x3 Fricke}
\bigl(t^2_3t_{12}+t^2_2t_{13}+t^2_1t_{23}
-(t_1t_2)^2-(t_2t_3)^2-(t_1t_3)^2-t_{321}'\bigr)&\notag\\
\times(t_{321}+ 
t^2_1+t^2_2+t_3^2&
-t_{12}-t_{13}-t_{23}) \\ 
=(t_{12}-t^2_1-t^2_2)(t_{13}-t^2_1-t^2_3)(t_{23}-t^2_2-&t^2_3). \notag
\end{align}
\end{rmk}

\subsection*{From $\bf3\times 3$ to $\bf2\times 2$} \ 

Now we will define a $\cB_3$-equivariant map from the space of 
triples of pseudo-reflections to the space of $\SL_2(\IC)$ triples. 
The main application of this here is just the observation that 
we will then obtain more exotic finite $\SL_2(\IC)$
braid group orbits from any triple of generators of a 
complex reflection group.

Suppose we are given the data 
$${\bt}:=(t_1,t_2,t_3,n_1,n_2,n_3,t_{12},t_{23},t_{13})$$
associated to a triple of pseudo-reflections.
(We extend the $\cB_3$-action to the set of such data---i.e. with square root
choices etc.---in the obvious way, permuting the $t_i$ and fixing the $n_i$.) 
Define a map $\varphi$ taking $\bf t$ to the 
$\SL_2(\IC)$ data $\bf m$ given by:


\begin{gather}
m_1:=\frac{t_1}{n_1}+\frac{n_1}{t_1},\qquad
m_2:=\frac{t_2}{n_1}+\frac{n_1}{t_2},\qquad
m_3:=\frac{t_3}{n_1}+\frac{n_1}{t_3}, \notag \\
m_{12}:= \frac{t_{12}}{t_1t_2},\qquad
m_{23}:= \frac{t_{23}}{t_2t_3},\qquad
m_{13}:= \frac{t_{13}}{t_1t_3},\label{eqn: varphi} \\
m_{321}:=\frac{n_2}{n_3}+\frac{n_3}{n_2}. \notag
\end{gather}

\begin{thm}\label{thm: B3 equivariance}
The map $\varphi$ is $\cB_3$-equivariant.
In particular finite $\cB_3$-orbits of $\SL_2(\IC)$ triples are obtained from 
triples of generators of three-dimensional complex reflection groups. 
\end{thm}
\pf
This may be proved by direct calculation (a less direct proof will be given in
section \ref{sn: imds}, 
along with a description of the origins of the above
formulae).
For example if we write 
$\bm'=\be_1(\varphi(\bt))$ and $\bm''=\varphi(\be_1(\bt))$, then
the tricky part is to see $m'_{12}=m''_{12}$. 
However it is straightforward to show that the expression obtained for  
$m'_{12}-m''_{12}$ (using the above formulae) has a factor of
$t_1t_2t_3-n_1n_2n_3$ in its numerator, which is zero due to 
\eqref{eq: sqrt of det}. (Similarly for $\be_2$.)
\epf

\begin{rmk}
It is possible to check directly (using Maple) that the map $\varphi$ is
well-defined; i.e. that $\varphi(\bt)$ satisfies the Fricke relation
\eqref{eq: fricke}
provided that $\bt$ satisfies both \eqref{eq: sqrt of det} and 
\eqref{eq: 3x3 Fricke}.
\end{rmk}

%
%

\subsection*{Painlev\'e parameters} \ 

The Painlev\'e VI equation that arises by performing
isomonodromic deformations of the rank $2$ Fuchsian system 
with monodromy data
$M_1,M_2,M_3,M_4$ (where $M_4M_3M_2M_1=1$) has parameters 
\begin{equation}\label{thal}
\al=(\theta_4-1)^2/2, \ \be=-\theta_1^2/2, \ 
\ga=\theta^2_3/2, \ \delta=(1-\theta_2^2)/2
\end{equation}
where the
$\theta_j$ $(j=1,2,3,4)$ are such that
$M_j \text{\ has eigenvalues } \exp(\pm\pi i\theta_j),$
i.e. $$\tr(M_j)=2\cos(\pi\theta_j).$$

Now suppose the $M_j$ arise under the map $\varphi$ from some data $\bt$
associated to a three-dimensional pseudo-reflection group.
We can then  relate the Painlev\'e parameters to the invariants of the 
pseudo-reflection group (cf. \cite{pecr} Lemma 3).
If we choose (for $j=1,2,3$) $\lambda_j,\mu_j$ such that 
\begin{equation}\label{eqn: bLambda}
t_j=\exp(\pi i \lambda_j),\qquad
n_j=\exp(\pi i \mu_j),\qquad \sum\lambda_i=\sum\mu_i
\end{equation}
then we have:
\begin{lem}\label{lem: pparams}
The Painlev\'e parameters corresponding to the data $\bt$ under the map
$\varphi$ are
\begin{equation}\label{thla}
\theta_i= \lambda_i-\mu_1\  (i=1,2,3),\quad\theta_4=\mu_3-\mu_2.
\end{equation}
\end{lem}
\pf
From the definition \eqref{eqn: varphi} 
of $\varphi$ we have $\tr(M_i)=m_i=2\cos\pi(\lambda_i-\mu_1)$ and 
$\tr(M_4)=m_{321}=2\cos\pi(\mu_3-\mu_2)$.
\epf

In particular if we are considering a complex reflection group $G$ 
that is generated by a triple of reflections,
then, by 5.4 of \cite{Shep-Todd},
we may choose a generating triple $(r_1,r_2,r_3)$ 
such that the $\mu_i$ are related to the exponents 
$x_1\leqslant x_2\leqslant x_3$ of the group $G$  
as follows:
$$ \mu_i= x_i/h,\quad (i=1,2,3)\qquad h:=x_3+1.$$ 
This result enables us to compile a table  
of parameters of the Painlev\'e
equations corresponding to the standard generators of those
irreducible three-dimensional
complex reflection groups which may be generated by three reflections 
(see \cite{pecr} Table 1, p.1019).



\subsection*{Example} \ 
Let us consider the reflections generating the Klein complex reflection group.
Explicitly the standard generators (from \cite{Shep-Todd} 10.1) are:

$$
r_1=\frac{1}{2}\left(\begin{matrix}
1&-1&-\overline{a} \\
-1&1&-\overline{a} \\
-a&-a&0
\end{matrix}\right),\qquad
r_2=\left(\begin{matrix}
1&0&0 \\
0&1&0 \\
0&0&-1
\end{matrix}\right),\qquad
r_3=\left(\begin{matrix}
1&0&0 \\
0&0&1 \\
0&1&0
\end{matrix}\right),
$$
where $a:=(1+i\sqrt{7})/2$.
These are order two complex reflections so $\det(r_i)=-1$ and
are ordered so that the exponents $\{3,5,13\}$
of the group appear in the eigenvalues
of the product $r_3r_2r_1$. Namely $r_3r_2r_1$ has eigenvalues
$\{
\exp(2\pi i\frac{3}{14}),
\exp(2\pi i\frac{5}{14}),
\exp(2\pi i\frac{13}{14}) \}$.
Also we compute:
$$\tr(r_1r_2)=1,\qquad
\tr(r_2r_3)=1,\qquad
\tr(r_1r_3)=0.$$
Thus if we set 
$\lambda_1=\lambda_2=\lambda_3=1/2, 
\mu_1=3/14, \mu_2=5/14, \mu_3=13/14$ so that 
$$t_1=t_2=t_3=i,\qquad 
n_1=\exp\frac{3\pi i}{14},\qquad
n_2=\exp\frac{5\pi i}{14},\qquad
n_3=\exp\frac{13\pi i}{14}$$
then the image of this data under $\varphi$ is
$$m_1=m_2=m_3=2\cos(2\pi/7),\quad m_{321}=2\cos(4\pi/7),\quad 
m_{12}=m_{23}=0,\quad m_{13}=1.$$
Clearly (cf. Lemma \ref{lem: pparams}) the parameters of the corresponding 
Painlev\'e equation are:
$$\theta_1=\theta_2=\theta_3=2/7,\quad\theta_4=4/7\qquad\text{and so }\quad(\al,\be,\ga,\de)=(9,-4,4,45)/98.$$

The corresponding braid group orbit is easy to calculate by hand;
Observe each of $m_1,m_2,m_3,m_{321}$ is fixed by $\cB_3$, and that,
since $4\cos(2\pi/7)\cos(4\pi/7)+4\cos^2(2\pi/7)=1$,
the formula for the action on the quadratic functions simplifies to
$$\be_1(m_{12},m_{23},m_{13})=
(1-m_{13}-m_{12}m_{23},m_{23},m_{12}),$$
$$\be_2(m_{12},m_{23},m_{13})=
(m_{12},m_{13},1-m_{23}-m_{13}m_{12}).$$
In this way we find the $\cB_3$ orbit has size seven, with values
\beq\label{bin no.s}
\begin {array}{ccc}
m_{12}&m_{23}&m_{13} \\
0&0&0\\
0&0&1\\
0&1&0\\
0&1&1\\
1&0&0\\
1&0&1\\
1&1&0.
\end{array}
\eeq
Upon restriction to the pure subgroup 
$\cP_3\subset \cB_3$, whose action is 
generated by $\be_1^2,\be_2^2$, we find the orbit still has
size seven. (More precisely, $\be_1^2,\be_2^2$ generate the free subgroup 
$\cF_2\subset \cP_3$, and $\cP_3$ is generated by $\cF_2$ and its centre,
which acts trivially.) 
Thus the corresponding solution of Painlev\'e VI has seven branches and,
via \eqref{bin no.s},
the branches of the solution are conveniently labelled by the
(binary) numbers from zero to six.
The generators of the pure braid group action
are represented by the following permutations of the seven
branches:
$$\be_1^2:\qquad(05)(14)(236)$$
$$\be_2^2:\qquad(03)(12)(465)$$
whose product also has two 2-cycles and a 3-cycle.
This should be the monodromy representation of the solution curve
as a cover of $\IP^1$ branched at $0,1,\infty$.
Using the Riemann--Hurwitz formula we thus 
see the solution curve has genus zero.
Also, for example, one can calculate the monodromy group of the cover;
the subgroup of the symmetric group
generated by these (even) permutations is as large as possible, namely $A_7$.
This gives a clear picture of the solution curve topologically as a cover of
the Riemann sphere. Our next aim (after explaining how the formulae
of this section were found) will be to find an explicit polynomial 
equation for the solution curve
and for the function $y$ on it solving Painlev\'e VI.


\begin{rmk}
Given the data corresponding to any of these branches one can easily solve the 
seven equations \eqref{eq: 7 2x2fns} to find a corresponding 
$\SL_2(\IC)$ triple.
For example for branch zero it is straightforward to find the triple:
\beq\label{eqn: Klein 2x2}
M_1=\left(\begin{matrix}
\phi&0 \\
0&\phi^{-1}
\end{matrix}\right),\qquad
M_2=\left(\begin{matrix}
w & x \\
-x &\overline{w}
\end{matrix}\right),\qquad
M_3=\left(\begin{matrix}
w&\mu x \\
-x/\mu&\overline{w}
\end{matrix}\right),
\eeq
where
$\phi=e^{2\pi i/7}$, $w=\frac{1+\phi^2}{\phi-\phi^3}$, 
$x=\sqrt{1-\vert w\vert^2}\in\IR_{>0}$ and
$\mu=(r+i\sqrt{4-r^2})/2$ where $r=\frac{(1+\phi^2)^2}{1-\phi^2}\in [0,1]
\subset\IR$.

\begin{lem}
The group generated by $M_1,M_2,M_3$ is an infinite subgroup of $\SU_2$.
\end{lem}
\pf
By construction the $M_i$ are in $\SU_2$ (somewhat surprisingly).
Since the group is nonabelian and contains elements of order seven we are
still to check it is not some large dihedral group.
Let $\varepsilon$ be an eigenvalue of $M_1^4M_2$,
so that 
$\varepsilon+\varepsilon^{-1}=\tau$ where 
$\tau=-(1+\phi^2)(\phi+\phi^4+\phi^6)$ and so
$\varepsilon=(\tau\pm\sqrt{\tau^2-4})/2$.
Thus $\varepsilon$ 
is some algebraic number of modulus one and we claim it is not a root
of unity.
To see this we take the product of $z-\varepsilon'$ over all the Galois conjugates of
$\varepsilon$ and find the minimal polynomial of $\varepsilon$ is
$$p(z)=z^6-3z^5-z^4-7z^3-z^2-3z+1.$$ 
Thus if $\varepsilon^N=1$ for some integer $N$ we would have that 
$p$ divides $z^N-1$ so all the roots of $p$ are roots of unity.
However $p(1)=-13$ so $p$ has a real root greater than $1$.
\epf

(Note also that triples in the same braid group orbit generate the same 
group.)

\end{rmk}

\subsection*{Some properties of ${\bf\varphi}$} \ 

Before moving on to the Fourier--Laplace transform we will describe some
properties of the map $\varphi$; we will show that the different choices of
square roots etc. give isomorphic  $\cP_3$ orbits and also examine
the fibres
of $\varphi$.
This motivates the definition of $\varphi$. 

First of all, the conjugacy class of a triple of pseudo-reflections only
determines the values of the functions
$t_1^2, t_2^2, t_3^2,
t_{12},
t_{23},
t_{13}$
and the symmetric functions in the $n_i^2$,
whereas the map $\varphi$ involves each $t_i,n_i$.
Thus we must make a choice of ordering of the eigenvalues $n_i^2$ of
$r_3r_2r_1$ and of the square roots $t_i,n_i$  such that $t_1t_2t_3=n_1n_2n_3$
in order to obtain the $2\times 2$ data. 
In general different choices lead to different $2\times 2$ data
(cf. Remark \ref{rmk sym. trips}). 
However the pure braid group orbits obtained via different choices are all
isomorphic:

\begin{lem}\label{lem: B3 equivariance under signs}
Let $\pi$ be a permutation of $\{1,2,3\}$ and choose signs 
$\varepsilon_i,\delta_i\in\{\pm1\}$ for $i=1,2,3$ such that
$\varepsilon_1\varepsilon_2\varepsilon_3=\delta_1\delta_2\delta_3$.
Consider the map $\sigma$ on the set of data
$${\bt}=(t_1,t_2,t_3,n_1,n_2,n_3,t_{12},t_{23},t_{13})$$
(satisfying \eqref{eq: sqrt of det} and \eqref{eq: 3x3 Fricke}) 
defined  by
$$\sigma(\bt):=
(\varepsilon_1t_1,\varepsilon_2t_2,\varepsilon_3t_3,
\delta_1n_{\pi(1)},\delta_2n_{\pi(2)},\delta_3n_{\pi(3)},
t_{12},t_{23},t_{13}).$$

Then map $\sigma$ commutes with the action of the pure braid group 
$\cP_3\subset
\cB_3$, and in particular
the $\cP_3$ orbits through $\varphi(\bt)$ and $\varphi(\sigma(\bt))$ are
isomorphic. 
\end{lem}
\pf
The pure braid group action is generated by $\be^2_1,\be_2^2$  and so fixes
the $t_i,n_i$ pointwise. From Lemma \ref{lem: 3x3 quad action} the action on
the functions $(t_{12},t_{23},t_{13})$ 
is independent of the sign and ordering choices.
\epf

We also wish to check that the different possible sign/ordering choices lead
to Painlev\'e VI equations with parameters which are equivalent under the
action of the affine $F_4$ Weyl group symmetries defined by
Okamoto \cite{OkaPVI}. 
To this end we lift the map $\sigma$ to act on the data 
${\bL}:=(\lambda_1,\lambda_2,\lambda_3,\mu_1,\mu_2,\mu_3)$ of 
\eqref{eqn: bLambda} as
$$\sigma(\bL)=
(\lambda_1+a_1,\lambda_2+a_2,\lambda_3+a_3,
\mu_{\pi (1)}+b_1,\mu_{\pi (2)}+b_2,\mu_{\pi (3)}+b_3)$$
where $a_i,b_i$ are integers such that $\sum a_i=\sum b_i$ and $\pi$ is a
permutation of $\{1,2,3\}$.
\begin{lem}\label{lem: Okamoto action on params}
The Painlev\'e VI parameters associated to 
${\bL}=(\lambda_1,\lambda_2,\lambda_3,\mu_1,\mu_2,\mu_3)$ 
in Lemma \ref{lem: pparams} are equivalent, under Okamoto's affine $F_4$ Weyl
group action, to those associated to $\sigma(\bL)$.
\end{lem}
\pf
Since the set of such $\sigma$'s forms a group it is sufficient  to check the
lemma on generators. The translations and the permutations may be dealt with
separately since the group is a semi-direct product.
First for the translations (fixing $\pi$ to be the
trivial permutation) this is straightforward; for example it is easy to
express the corresponding translations of the $\theta$'s in terms of the
translations of \cite{NY-so8} (34).
Finally the permutation just swapping $\mu_1$ and $\mu_2$ 
is obtained from the
transformation $s_1s_2s_1$ (of \cite{NY-so8}) and 
that swapping $\mu_1$ and $\mu_3$ is obtained from the
transformation $(s_0s_1s_3s_4)s_2(s_0s_1s_3s_4)$.
\epf

\begin{rmk}
Perhaps it is helpful to recall that there are several symmetry groups of PVI
considered in the literature, amongst which we have
$$\text{affine } D_4\  <\  \text{extended affine } D_4\  <\  
\text{affine } F_4$$
the first two of which are for example considered in \cite{NY-so8}.
In brief\,\footnote{I am grateful to M. Noumi for clarifying this to me.} 
the first two differ by the Klein four-group and do not involve
changing the time parameter $t$, whereas the full affine $F_4$ action of
Okamoto involves changing $t$ (by automorphisms of $\IP^1$ permuting
$0,1,\infty$). In fact only the extended affine $D_4$ symmetries were used
above, although the full $F_4$ action will be considered in section
\ref{sn: ineq}. 
\end{rmk}

Next we will examine the fibres of $\varphi$. 
From our rough dimension counts we see these fibres should be one dimensional.
The continuous part of the fibres arises as follows.
Define an action of $\IC^*$ on the pseudo-reflection
data $\{\bt\}$ by declaring $h\in\IC^*$ to act as
\beq\label{eqn: C* action defn downstairs}
t_i\mapsto ht_i, \qquad 
n_i\mapsto hn_i,\qquad
t_{ij}\mapsto h^2t_{ij}.
\eeq
Observe that this does indeed act within the fibres of $\varphi$, i.e.
that $\varphi(h\bt)=\varphi(\bt)$ for any $h\in\IC^*$.
Moreover a simple direct calculation shows this $\IC^*$ action commutes with
the $\cB_3$ action on $\{\bt\}$.
(The simplicity of this action is deceptive since we carefully 
chose the functions $t_{i},t_{ij}$.) 

Thus, for example, we can always use this action to move to the
($\cB_3$-invariant) subset
of the pseudo-reflection data having $n_1=1$.
\begin{lem}
The map $\varphi$ is surjective, and the restriction of $\varphi$ to the
subset of the pseudo-reflection data $\{\bt\}$ having $n_1=1$ is a finite map.
\end{lem}
\pf
Given arbitrary $\SL_2(\IC)$ data $\bm$ we just try to solve for $\bt$ (having
first set $n_1=1$). One finds a solution always exists and there are five sign
choices, so a generic fibre has $32$ points.\epf
\begin{rmk}
One  may check algebraically that if $\bt$ has $n_1=1$,
satisfies $n_2n_3=t_1t_2t_3$ and is such that $\varphi(\bt)$
satisfies the Fricke relation \eqref{eq: fricke}
then $\bt$ satisfies the $3\times 3$ analogue \eqref{eq: 3x3 Fricke}
of the Fricke relation.
\end{rmk}

%
%
%
%
%
%
%
%
%

Moving to the subset of the data on which $n_1=1$ implies we are forcing
$1$ into the spectrum of the product $r_3r_2r_1$.
This implies that
the representation (of the free group on three letters) defined
by $(r_3,r_2,r_1)$ is reducible:
This is clear if $r_i=1+e_i\otimes\al_i$ for some 
$e_i$ which are not a basis of $\IC^3$ (since the span of the $e_i$ is an
invariant subspace). Otherwise we have:
\begin{lem}[cf. \cite{Carter} 10.5.6, \cite{Coleman} 3.7]\label{lem:redtrip}
If $r_i:=1+e_i\otimes\al_i$ for a basis $e_1,e_2,e_3$ of $V=\IC^3$ and
$v\in V$ satisfies $r_3r_2r_1v=v$ then
$r_iv=v$ for $i=1,2,3$.
\end{lem}
\pf
If $r_3r_2r_1v=v$ then $$r_2r_1v-v=r^{-1}_3v-v,$$ 
the lefthand side of which is
a linear combination of  $e_2,e_1$, and the righthand side is a multiple of 
$e_3$, since $r_3^{-1}=1-e_3\otimes\al_3/t_3^2$.  
Thus both sides vanish so $r_3v=v$ and $r_2r_1v=v$. 
Then similarly we see both sides $r_1v-v=r_2^{-1}v-v$ vanish.
\epf

Thus we can use the
$\IC^*$ action to move to a reducible triple.
The map $\varphi$ is defined simply by first moving to a reducible triple
by setting $h:=n_1^{-1}$ so
that we can write the $r_i$ in block upper triangular form.
In general there will then 
be a size two and a size one block (with entry $1$) on the diagonal.
We then define $\wh M_i\in \GL_2(\IC)$ to be the size two block,
which will have eigenvalues
$\{1,(t_i/n_1)^2\}$.
Hence defining $M_i=n_1\wh M_i/t_i$ yields an $\SL_2(\IC)$ triple.
Computing the various traces then  gives the stated formulae for the map
$\varphi$. 
(In more invariant language we take the projection to $\SL_2(\IC)$ of the
rank two part of the `semisimplification' of the reducible representation.)

Thus we have motivated the map $\varphi$ in terms of the $\IC^*$ action.
In section \ref{sn: imds} we will motivate this action as the image under the
Fourier--Laplace transform 
of a simple scalar shift. 

\subsection*{Orthogonal reflection groups} \ 

Let us check that the unipotent $2\times 2$ monodromy data and the orthogonal
three-dimensional reflection groups considered by Dubrovin--Mazzocco 
in \cite{DubMaz00} are
related by the map $\varphi$ defined in \eqref{eqn: varphi}.

The monodromy data in \cite{DubMaz00} is parameterised by four numbers 
$(x_1,x_2,x_3,\mu)$ related by the condition 
$$m=2\cos(2\pi \mu) \qquad\text{ where } m:=2+x_1x_2x_3-(x_1^2+x_2^2+x_3^2),$$
which is equivalent to \cite{DubMaz00} (1.21). The degenerate cases 
$m=\pm 2$ are ruled out.
Although the main interest is in real orthogonal
reflection groups the formulae here make
sense for complex values of the parameters; effectively we are restricting
$\varphi$ to a complex three-dimensional slice.

Without loss of generality one may assume $x_1\ne 0$ and then 
the $2\times 2$ monodromy data is given by the triple 
(\cite{DubMaz00} (1.20)):
\beq\label{eqn: unip}
M_1=
\left(\begin{matrix}
1 & -x_1  \\
0 & 1
\end{matrix}\right),\qquad
M_2=\left(\begin{matrix}
1 & 0  \\
x_1 & 1
\end{matrix}\right),\qquad
M_3=\left(\begin{matrix}
1+x_2x_3/x_1 & -x_2^2/x_1  \\
x^2_3/x_1 & 1-x_2x_3/x_1
\end{matrix}\right).
\eeq
Note that, generically,  each of these matrices is conjugate to 
$\left(\begin{smallmatrix}
1 & 1  \\
0 & 1
\end{smallmatrix}\right)$, and in all cases $m_j:=\tr(M_j)=2$.
Also, straightforward computations give that
$$
m_{12}=2-x_1^2,\qquad
m_{23}=2-x_2^2,\qquad
m_{13}=2-x_3^2,\qquad
m_{321}=m.$$
(By nondegeneracy there is always some index $j$ such that $x_j\ne 0$ and one
may obtain the same values of the invariant functions from an analogous 
$2\times 2$ triple.)

The corresponding 
$3\times 3$ reflections considered in \cite{DubMaz00} are
(\cite{DubMaz00} (1.51)):
$$
r_1=\left(\begin{matrix}
-1 & -x_1 & -x_3  \\
0&1&0\\
0&0&1
\end{matrix}\right),\qquad
r_2=\left(\begin{matrix}
1&0&0\\
-x_1 & -1 & -x_2  \\
0&0&1
\end{matrix}\right),\qquad
r_3=\left(\begin{matrix}
1&0&0\\
0&1&0\\
-x_3 & -x_2 & -1  
\end{matrix}\right)
$$
which preserve the nondegenerate symmetric bilinear form
given by the matrix 
$$\left(\begin{matrix}
2 & x_1 & x_3 \\
x_1 &  2 & x_2 \\
x_3 & x_2 & 2 
\end{matrix}\right).$$
Immediate computation then gives that
$$
t_1^2=t_2^2=t_3^2=-1,$$
$$
t_{12}=x_1^2-2,\qquad
t_{23}=x_2^2-2,\qquad
t_{13}=x_3^2-2$$
$$t_{321}=-m-1=-t'_{321}.$$
The characteristic polynomial of $r_3r_2r_1$ is
$$(\lambda+1)(\lambda^2+m\lambda+1)$$
which has roots 
$$n_1^2=-1,\qquad n_2^2= -\exp(2 \pi i \mu),\qquad n_3^2= -\exp(-2 \pi i\mu).$$

Now we claim that if we choose the square roots appropriately then the map
$\varphi$ of \eqref{eqn: varphi} takes this $3\times 3$ data onto the
unipotent $2\times 2$ data above.
Indeed setting $n_1=t_1=t_2=t_3=i$ clearly gives the correct values
$m_j=i/i+i/i=2$ and $m_{jk}=t_{jk}/i^2=-t_{jk}$.
Also if we set $n_2=i\exp(\pi i \mu)$ and $n_3=i\exp(-\pi i \mu)$ then
$n_2/n_3=\exp(2\pi i \mu)$ and so
$m_{321}=n_2/n_3+n_3/n_2=2\cos 2\pi\mu=m$ as required.
Thus the map $\varphi$ does indeed extend the above correspondence used by
Dubrovin--Mazzocco. 

\begin{rmk} \label{rmk sym. trips}
If instead we choose to order the eigenvalues of $r_3r_2r_1$ as
$$n_1^2=-\exp(2 \pi i \mu) ,\qquad n_2^2=-1,\qquad n_3^2= -\exp(-2 \pi i\mu)$$
then we claim that, 
with appropriate square root choices, the corresponding $2\times 2$ data
(under $\varphi$) 
has the remarkable property that the four local monodromies
$M_1,M_2,M_3, M_3M_2M_1$ all lie in the same conjugacy class:
Namely if we choose $n_2=t_1=t_2=t_3=i$, let $n_1$ be any square root of 
$-\exp(2 \pi i \mu)$ and define $n_3:=1/n_1$
then we find
$$m_1=m_2=m_3=m_{321}=i/n_1+n_1/i=\pm2\cos(\pi \mu).$$
Thus if this common value is not $\pm 2$, the corresponding $\SL_2(\IC)$
matrices are regular semisimple and we have established the claim.
However if this value is $\pm 2$, it follows that $m=2$
contradicting the nondegeneracy assumption.

Such `symmetric' $\SL_2(\IC)$ triples have been studied in this context by
Hitchin (cf. \cite{Hit-Poncelet, Hit-Octa}), 
and that they arise from real reflection groups was known to
Dubrovin--Mazzocco (cf. \cite{DubMaz00} Remark 0.2). 
It is interesting to note that, from triples of generating reflections of the
real three-dimensional tetrahedral, octahedral, and icosahedral reflection
groups, one
obtains in this way 
triples of generators of finite subgroups of $\SL_2(\IC)$. 
However, rather bizarrely, the tetrahedral and octahedral 
groups are
swapped in the process: the tetrahedral reflection group maps to the binary
octahedral subgroup of $\SL_2(\IC)$ and 
vice-versa; the  octahedral reflection group maps to the binary
tetrahedral subgroup of $\SL_2(\IC)$. The three inequivalent triples of
generating reflections of the icosahedral reflection group do all map to
triples of generators of the binary icosahedral group though. 
\end{rmk}

\end{section}

\begin{section}{Isomonodromic deformations}\label{sn: imds}

The main aim of this section is
to see how the $\IC^*$ action of \eqref{eqn: C* action defn downstairs}
on the invariants of  the pseudo-reflection data
arises, and, on the other side of the Riemann--Hilbert
correspondence, to describe the corresponding $\IC$ action
on the rank three Fuchsian systems.

This is the key ingredient needed to motivate the map $\varphi$ of
\eqref{eqn: varphi}, as described after Lemma \ref{lem:redtrip} above.
As a corollary we will see why $\varphi$ is $\cB_3$-equivariant.

We will work in a somewhat more general context in this section than the rest
of the paper;
the reader interested mainly in the construction of algebraic solutions
to PVI could skip straight to the next section.

Apart from the desire to explain how the procedure of section 
\ref{sn: braid orbits} was found, the motivation for this section is to
enable us to see (in section \ref{sn: reconstruction})
how one may work back from an explicit solution to PVI to an explicit 
rank three system of differential equations. 
This will give a mechanism for constructing new non-rigid systems of
differential equations with finite monodromy group.
(Except for this 
the proofs given in the other sections are independent of the results of 
this section.) 

Let us begin with some generalities on isomonodromic deformations of Fuchsian
systems and the Schlesinger equations.
Let $V=\IC^n$ and suppose we have 
matrices 
$B_1,\ldots,B_{m-1}\in\End(V)$ and distinct points
$a_1,\ldots, a_{m-1}\in\IC$.
Then consider
the following meromorphic
connection on the trivial rank $n$ holomorphic vector bundle over the
Riemann sphere:
\begin{equation} \label{nabla}
\nabla := d - \left(B_1\frac{dz}{z-a_1}+\cdots+B_{m-1}\frac{dz}{z-a_{m-1}}
\right).
\end{equation}
This has a simple pole at each $a_i$ and at infinity.
Write 
$$B_m=B_\infty:=-(B_1+\cdots+B_{m-1})$$
for the residue matrix at infinity.
Thus, on removing disjoint open discs $D_1,\ldots,D_m$ 
from around the poles and
restricting $\nabla$ to the $m$-holed sphere
$${{S}}:=\IP^1\setminus(D_1\cup\cdots\cup D_m),$$
we obtain a (nonsingular)
holomorphic connection.
In particular it is flat and so, taking its monodromy, a
representation of the fundamental group of the $m$-holed sphere is
obtained.
This procedure defines a holomorphic map, which we will call the
monodromy map or Riemann--Hilbert map, from the set of such connection 
coefficients to the set of complex
fundamental group representations:
$$\bigl\{ (B_1,\ldots,B_m) \ \bigl\vert \ \text{$\sum B_i= 0$} \bigr\}
\ \mapright{\text{RH}}\ 
\bigl\{ (M_1,\ldots,M_m) \ \bigl\vert \ M_m\cdots M_1 = 1 \bigr\}
$$
where appropriate loops generating the fundamental group
of ${{S}}$ have been chosen 
and the matrix $M_i\in G:=\GL_n(\IC)$ is the
automorphism obtained by parallel translating a basis of solutions
around the $i$th loop.

The Schlesinger equations are the equations for isomonodromic deformations 
of the connection \eqref{nabla}. Suppose we move the pole positions 
$a_1,\ldots,a_{m-1}$. Then we wish to vary the coefficients $B_i$, as
functions of the pole positions, such that the monodromy data 
$(M_1,\ldots,M_m)$ only changes by diagonal conjugation by $G$.
This is the case if the $B_i$ 
vary according to Schlesinger's equations:
\beq\label{eq: schles}
\frac{\partial B_i}{\partial a_j}=\frac{[B_i,B_j]}{a_i-a_j}\qquad 
\text{if } i\ne j, \quad \text{ and} \qquad 
\frac{\partial B_i}{\partial a_i}=
-\sum_{j\ne i,m}\frac{[B_i,B_j]}{a_i-a_j}
\eeq
where $i=1,\ldots,m-1$.
Observe that these equations imply that the
residue at infinity $B_m$ is held constant.
Also note that 
the Schlesinger equations are equivalent to the flatness of the connection
\begin{equation} \label{eqn: bignabla}
d - \left(B_1\frac{dz-da_1}{z-a_1}+
\cdots+B_{m-1}\frac{dz-da_{m-1}}{z-a_{m-1}}
\right).
\end{equation}

In terms of differential forms Schlesinger's equations may be rewritten
as 
\beq
dB_i = - \sum_{j\ne i,m}[B_i,B_j]d_{ij}
\eeq
where $d$ is the exterior derivative on $\{a_i\}$ and
$d_{ij}:=d\log(a_i-a_j)=(da_i-da_j)/(a_i-a_j)$.
In turn it will be convenient to rewrite this as 
\beq\label{eq 3rd schles}
dB_i=[L_i,B_i]\qquad\text{where }\quad L_i:=\sum_{j\ne i,m}{B_j d_{ij}}.
\eeq

Note that if we have a local solution of Schlesinger's equations and we
construct the $L_i$ from the formula \eqref{eq 3rd schles} then
firstly we have that $\nabla_i:=d-L_i$ is a flat connection and secondly that
$B_i$ is a horizontal section of $\nabla_i$ (in the adjoint representation).

Now let us specialise to the case where the dimension 
$n$ equals the number $m-1$
of finite singularities, and where each of the finite residues 
$B_1,\ldots, B_n$ is a rank one matrix.
Thus 
$$B_i=f_i\otimes\be_i\qquad\text{for some}\quad f_i\in V,\  \be_i\in V^*.$$ 

Then we may lift the Schlesinger equations from the space of residues $B_i$ to
the space of $f_i$'s and $\be_i$'s.
Namely, suppose we have a local solution of the  Schlesinger equations 
on some polydisc $D$. 
Then we can write $B_i=f_i\otimes\be_i$ for $i=1,\ldots,n$ at some base-point
and evolve $f_i,\be_i$ over $D$, as solutions to
the equations:
\beq\label{eq: lifted schles.}
df_i=L_if_i\qquad d\be_i=-\be_iL_i
\eeq 
where the $L_i$ are defined in terms of the given $B_j$ solving the Schlesinger
equations.
Then one finds immediately that the $f_i\otimes \be_i$ solve 
\eqref{eq 3rd schles}, and so $f_i\otimes \be_i=B_i$ throughout $D$ (since
they agree at the basepoint and solve Schlesinger's equations).
Alternatively one can view \eqref{eq: lifted schles.} 
as a coupled system of nonlinear equations for $\{f_i,\be_i\}$, by defining 
$L_i$ in terms of the $B_j:=f_j\otimes \be_j$ as in \eqref{eq 3rd schles}.
We will refer to these as the {\em lifted equations} 
(they were introduced in \cite{JMMS} and further studied in \cite{Harn94}).
The above considerations show:
\begin{prop}
Any solution of the Schlesinger equations \eqref{eq: schles} may be lifted to
a solution of the lifted equations \eqref{eq: lifted schles.}
by only solving linear equations.
Conversely any solution of the lifted equations projects to a
solution of \eqref{eq: schles} by setting $B_i= f_i\otimes\be_i$.
\end{prop}

Now we wish to define an action of $\IC$ which will be
the additive analogue of the 
$\IC^*$ action of \eqref{eqn: C* action defn downstairs}.

Suppose we have a local solution
$\{B_1({\bf a}),\ldots,B_n({\bf a})\}$ 
of Schlesinger's equations on some polydisc $D$,
where ${\bf a}=(a_1,\ldots,a_n)$,
such that the images of the $B_i$ are linearly independent (i.e. for any
choice of $f_i,\be_i$ such that $B_i=f_i\otimes\be_i$, the $f_i$ make up a
basis of $V$).
Then we can define the following action of the complex numbers on the set of
such solutions:
\begin{prop}\label{prop: additive translation}
For any complex number $\lambda\in\IC$ the matrices
$$\wt B_i:= B_i+\lambda f_i\otimes \wh f_i$$
constitute another solution to Schlesinger's equations on $D$, where
$\wh f_1,\ldots,\wh f_n\in V^*$ are the dual basis defined by 
$\wh f_i(f_j)=\delta_{ij}$.
\end{prop}
\pf
First note that this is well-defined since the projectors
$f_i\otimes \wh f_i$ are independent of the choice of $f_i$'s.
Then lift the $B_i$ arbitrarily 
to a solution $\{f_i({\bf a}),\be_i({\bf a})\}$ of 
\eqref{eq: lifted schles.} over $D$.
Straightforward computations then give that $\wh f_i$ satisfies
$d\wh f_i=-\sum_{j\ne i}\be_i(f_j)\wh f_j d_{ij}$
and using this one easily confirms 
$d\wt B_i=[\wt L_i,\wt B_i]$ 
where
$\wt L_i=  L_i+\lambda\sum_{j\ne i}f_j\otimes \wh f_j d_{ij}.$
\epf

One may arrive at this action 
as follows.
Given a local solution $\{f_i({\bf a}),\be_i({\bf a})\}$ 
of the lifted equations
one may check that the matrix $B\in \End(V)$
defined by $(B)_{ij}=\be_i(f_j)$
satisfies the nonlinear differential equation
\beq\label{eq: dual eqn}
dB= \left[B, \ad^{-1}_{A_0}([dA_0,B])\right]
\eeq
where $A_0:=\diag(a_1,\ldots,a_n)$.
(Note that $\ad_{A_0}:\End(V)\to \End(V)$ is invertible when restricted to the
matrices with zero diagonal part and that $[dA_0,B]$ has zero diagonal part.)
This is the `dual' equation to the Schlesinger equations in the 
present context (in
the sense of Harnad \cite{Harn94}) 
and arises as the equation for isomonodromic deformations of the 
irregular connection
\beq\label{eq: irr conn}
d-\left(\frac{A_0}{w^2}+\frac{B}{w}\right)dw,
\eeq
which, after an appropriate coordinate change, appears as 
the (twisted) Fourier--Laplace transform of the original Fuchsian system 
(cf. \cite{BJL81} and references therein).
Equation \eqref{eq: dual eqn} 
appears in the theory of Frobenius manifolds \cite{Dub95long}
for skew-symmetric $B$ and is related to quantum Weyl groups \cite{bafi}.
Note that equation \eqref{eq: dual eqn} 
is equivalent to the Schlesinger equations in that its solutions may  also be
lifted to solutions of \eqref{eq: lifted schles.} by only solving the linear
equations 
$$
df_i=   \sum_{j\ne i}(B)_{ji}f_jd_{ij}\qquad 
d\be_i=-\sum_{j\ne i}(B)_{ij}\be_jd_{ij}$$
where $B({\bf a})$ solves \eqref{eq: dual eqn}.

Now from the form of \eqref{eq: dual eqn} it is transparent that 
replacing $B$ by $B+\lambda$ maps solutions to solutions (where 
$\lambda\in \IC$ is
constant). 
(Observe this corresponds to tensoring the irregular connection 
\eqref{eq: dual eqn} by the meromorphic 
connection $d-\lambda\frac{dw}{w}$ on the trivial line bundle.)
If $B$ is translated in this way, then (provided the $f_i$ are a basis) 
we can see how to change the
corresponding Schlesinger solutions as follows.
First note:
\begin{lem}\label{lem: residue conjugacy}
Suppose $f_1,\ldots,f_n\in V$ is an arbitrary basis, 
$\be_i\in V^*$ is arbitrary, $B_i=f_i\otimes\be_i\in\End(V)$
and $(B)_{ij}=\be_i(f_j)$.
Then 
$$(B_1,\ldots,B_n)\text{ is conjugate to }
(E_1B,\ldots,E_nB)$$
where $E_i\in\End(V)$ has $(i,i)$ entry $1$ and is otherwise zero. 
\end{lem}
\pf
Define $g\in \GL(V)$ to have $i$th column $f_i$.
Then observe that $g^{-1}B_ig=E_iB$.
\epf

Thus replacing $B$ by $B+\lambda$ changes $B_i=gE_iBg^{-1}$ to 
$B_i+\lambda gE_ig^{-1}=B_i+\lambda f_i\otimes\wh f_i$
and so we deduce the action of Proposition \ref{prop: additive translation}.

Now the basic idea to reduce the rank of the systems 
by one is to choose $\lambda$ to be an
eigenvalue of $B_\infty$. Then the residue at infinity 
$\wt B_\infty$
of the resulting system has a nontrivial kernel. 
This is because
Lemma \ref{lem: residue conjugacy} implies
\begin{equation}\label{stmt: conjugate residues}
B_\infty=-(B_1+\cdots +B_n)\quad\text{is conjugate to}\quad 
-(E_1B+\cdots+E_nB)=-B.
\end{equation}
Thus translating $B$ by $\lambda$ implies that 
$\wt B_\infty$ is conjugate to $B_\infty-\lambda$, which will have a zero
eigenvalue. 
Say $\wt B_\infty v=0$ for some vector $v$. Then the fact that 
the $f_i$ are a basis implies that $v$ is in the kernel of {\em all}
the residues $\wt B_i$---the 
resulting system is reducible and we can pass to the
corresponding $(n-1)\times (n-1)$ quotient system.

The next step is to find the action on monodromy data 
corresponding to the $\IC$ action above.
Suppose $B_j=f_j\otimes\be_j$ for $j=1,\ldots,n$ and each 
$$\lambda_j:=\tr(B_j)=\be_j(f_j)$$
is not an integer.
Then one knows that the monodromy matrix $M_j$ 
around $a_j$ is conjugate to $\exp(2\pi i B_j)$ and so is a 
(diagonalisable) pseudo-reflection.
We will write
$$r_j=M_j=1+e_j\otimes \al_j\qquad\text{where $e_j\in V,\al_j\in V^*$}$$
for this pseudo-reflection.
Clearly the non-identity eigenvalue of $r_j$ is $\exp(2\pi i \lambda_j)$
so setting $t_j=\exp(\pi i \lambda_j)$ (as in \eqref{eqn: bLambda}) implies
$t_j^2=\det(r_j)$ agreeing with the definition \eqref{eqn: t defn}.

From \eqref{stmt: conjugate residues} we deduce that if
$B$ is translated by $\lambda$
then the monodromy around a large positive loop is just scaled:
\beq\label{eq: monod scaled}
\wt r_n\cdots \wt r_2\wt r_1\quad\text{is conjugate to}\quad 
r_n\cdots r_2r_1h^2
\eeq
where $$h:=\exp(\pi i \lambda)\in\IC^*,$$
at least if $B_\infty$ is sufficiently generic (no distinct eigenvalues
differing by integers).
(Here $\wt r_1,\ldots,\wt r_n$ are the monodromy data of the connection
obtained by replacing each $B_i$ by $\wt B_i$ in \eqref{nabla}.)

In brief the additive action was determined by the fact that sum
$\sum B_i$ was just translated by $\lambda$ (assuming the $f_i$  make up a
basis, which is held fixed).
We will see below that 
the multiplicative action (i.e. the action on monodromy data) is determined
by the fact that the product $r_n\cdots r_1$ is just scaled by $h^2$ 
(assuming the $e_i$ make up a basis, which is held fixed).

First let us recall a basic algebraic fact about pseudo-reflections.

Suppose $e_1,\ldots, e_n$ are a basis of $V$ and 
$\al_1,\ldots,\al_n\in V^*$ are such that
$r_i:=1+e_i\otimes \al_i\in \GL(V)$, i.e. $1+\al_i(e_i)\ne 0$.
Define two $n\times n$ matrices $t^2,u$ by
$$
t^2:=\diag(1+\al_1(e_1),\ldots,1+\al_n(e_n)),\qquad 
(u)_{ij}=\al_i(e_j).
$$
(We do not need to choose a square root $t$ of $t^2$ at this stage, but it is
convenient to keep the notation consistent with other sections of the paper.)
Then let $u_+,u_-$ be the two unipotent matrices 
determined by the equation
\beq \label{eq u,up,um}
t^2u_+-u_-=u              
\eeq
where $u_+$ is upper triangular with ones on the diagonal and 
$u_-$ is lower triangular with ones on the diagonal.

\begin{thm}
[Killing \cite{killing2}, 
Coxeter \cite{coxeter51}] \label{thm: killing}
The matrix representing the product $r_n\cdots r_1$ (in the $e_i$ basis)
is in the big-cell of $\GL_n(\IC)$, and so may be written uniquely as the
product of a lower triangular, a diagonal and 
an upper triangular matrix.
Moreover this factorisation is given explicitly by $u_-^{-1}t^2u_+$:
\beq\label{eq: killing's eqn}
r_n\cdots r_2r_1=u_-^{-1}t^2u_+.
\eeq
\end{thm}

\begin{rmk}
The history of this 
result is discussed by Coleman \cite{Coleman} (cf. Corollary 3.4).
Coxeter
proves this for genuine reflections---i.e. coming from 
a symmetric bilinear form in \cite{coxeter51}.
The starting point of this paper was the simple observation 
that Coxeter's argument may be extended to the pseudo-reflection case.
Dubrovin had used Coxeter's version in relation to Frobenius manifolds 
(cf. \cite{DubPT}) and the author was interested in extending Dubrovin's
picture to the general case (cf. \cite{smapg}).
Despite asking various complex reflection group experts 
the author only  found Coleman's paper (and hence the link to Killing) 
since \cite{Coleman} is in the same volume as a 
well-known paper of Beukers--Heckman.
\end{rmk}

It is worth clarifying the fact that generically the matrix $u$ determines
$(r_1,\ldots,r_n)$ up to conjugacy:
\begin{lem}
If $\det(u)\ne 0$ then there 
is a matrix $g\in \GL(V)$
such that, for $i=1,\ldots,n$ we have 
$$r_i= g(1+e_i^o\otimes\ga_i)g^{-1}$$
where
$\ga_i\in V^*$ 
is the $i$th row of the matrix $u$ 
and $e_i^o$ is the standard basis of $V$.
\end{lem}
\pf
By definition of $u$, if 
$\det(u)\ne 0$ then the $e_i$ are a basis of $V$.
Then the result follows since we know the action on a basis:
$r_i(e_j)=e_j+u_{ij}e_i$.
\epf

Note that if we define $u_\pm,t^2$ by the equation 
$u=t^2u_+-u_-$ then the condition that $\det(u)\ne 0$ is equivalent to saying 
$1$ is not an eigenvalue of $u_-^{-1}t^2u_+$, since
$\det(u)=\det(t^2u_+-u_-)=\det(u_-^{-1}t^2u_+-1)$.

Thus generically the $n$-tuple $(r_1,\ldots,r_n)$ is determined
up to overall conjugacy by the matrix $u$, and in turn by the product
$r_n\cdots r_1$, by Theorem \ref{thm: killing}.
From \eqref{eq: monod scaled} the obvious guess is therefore that
that $\wt r_n\cdots \wt r_1= u_-^{-1}t^2u_+h^2=u_-^{-1}h^2t^2u_+$ 
so that $\wt \al_i(\wt e_j)= (h^2t^2u_+-u_-)_{ij}$, 
which should determine
$(\wt r_1,\ldots,\wt r_n)$ up to overall conjugacy. 
The following theorem says that this is indeed the case, at least generically.
Suppose each $\lambda_i$ and each eigenvalue of $\sum B_i$ is not an integer
(and that the same holds after translation by $\lambda$).
Then we have:

\begin{thm}[Balser--Jurkat--Lutz \cite{BJL81}]\label{thm: FL on monod}
Let $\wt u$ be the matrix
$$\wt u = h^2t^2u_+-u_-$$ 
where $u_\pm,t^2,h$ are as defined above.
Then there is a basis $\wt e_i$ of $V=\IC^n$ 
and $\wt\al_1,\ldots\wt\al_n\in V^*$
such that $\wt r_i=1+\wt e_i\otimes\wt \al_i$ and
$(\wt u)_{ij}=\wt \al_i(\wt e_j)$ for all $i,j$.
\end{thm}

\begin{rmk}
This is not written down in precisely this way in \cite{BJL81}
so we will describe how to extract it in Appendix 
\ref{apx: bjl extraction}. 
The key point is that the matrices $u_\pm$ are essentially
the Stokes matrices of the
irregular connection \eqref{eq: irr conn} and are easily seen to be
preserved under the scalar shift.
Then one computes bases of solutions of the Fuchsian connection 
as Laplace transforms of standard bases of solutions of \eqref{eq: irr conn}
and this enables the Stokes matrices to be related to the 
pseudo-reflection data
$u$ as in equation \eqref{eq u,up,um}.
The observation that this implies the Fuchsian monodromy data and the Stokes
data are 
then related by the beautiful equation \eqref{eq: killing's eqn} in
Theorem  \ref{thm: killing} does not seem to appear in \cite{BJL81}. 
In summary we see that equation \eqref{eq: killing's eqn}
is the manifestation of the Fourier--Laplace transformation 
on monodromy data, relating the monodromy data of the Fuchsian connection to 
the monodromy/Stokes data of the corresponding irregular connection.
\end{rmk}

In other words: in general 
the matrix $u$ determines $(r_1,\ldots,r_n)$ up to overall
conjugation and Theorem \ref{thm: FL on monod} explains how the matrix $u$ 
varies: the lower triangular part is fixed, the upper triangular part is
scaled by $h^2$, and the diagonal part $t^2-1$ becomes $h^2t^2-1$.
Let us make this more explicit in the $n=3$ case.
We start with a connection 
\beq\label{eq: rk3 fuchsian} 
d-\sum_1^3\frac{B_i}{z-a_i}dz
\eeq
where, up to overall conjugation:
$$
B_1=\left(\begin{matrix}
\lambda_1& b_{12} & b_{13} \\
0&0&0\\
0&0&0
\end{matrix}\right),
\qquad
B_2=\left(\begin{matrix}
0&0&0 \\
b_{21}&\lambda_3&b_{23} \\
0&0&0
\end{matrix}\right),\qquad
B_3=\left(\begin{matrix}
0&0&0 \\
0&0&0 \\
b_{31}&b_{32}&\lambda_3
\end{matrix}\right)
$$
for some numbers $b_{ij}$ with $i\ne j$.
Then we take the monodromy data of this and obtain pseudo-reflections
$r_1,r_2,r_3$ which, up to overall conjugation, are of the form
$$
r_1=\left(\begin{matrix}
t^2_1&u_{12}&u_{13}\\
0&1&0 \\
0&0&1
\end{matrix}\right),\qquad
r_2=\left(\begin{matrix}
1&0&0 \\
u_{21}&t_2^2&u_{23} \\
0&0&1
\end{matrix}\right),\qquad
r_3=\left(\begin{matrix}
1&0&0 \\
0&1&0 \\
u_{31}&u_{32}&t^2_3
\end{matrix}\right)
$$
where
$t_j=\exp(\pi i \lambda_j)$.
Then we replace $\lambda_i$ by $\lambda_i+\lambda$ in 
\eqref{eq: rk3 fuchsian} for each $i$, and Theorem \ref{thm: FL on monod}
says that the monodromy of the resulting connection is conjugate to
\beq\label{eq: translated monod}
\wt r_1=\left(\begin{matrix}
h^2t_1^2 &h^2u_{12}&h^2u_{13}\\
0&1&0 \\
0&0&1
\end{matrix}\right),\ 
\wt r_2=\left(\begin{matrix}
1&0&0 \\
u_{21}&h^2t_2^2&h^2u_{23} \\
0&0&1
\end{matrix}\right),\ 
\wt r_3=\left(\begin{matrix}
1&0&0 \\
0&1&0 \\
u_{31}&u_{32}&h^2t^2_3
\end{matrix}\right)
\eeq
where $h:=\exp(\pi i \lambda)$.


Taking the various traces yields the fact that the invariant functions of the
monodromy matrices are related as:
$$
{\wt t_i}^2=h^2t_i^2,\qquad 
\wt t_{ij}=h^2t_{ij},\qquad
\wt t_{321}=h^2t_{321},\qquad 
\wt t'_{321}=h^4t'_{321}.
$$
These equations hold for any $\lambda$
since the invariants are analytic functions of the coefficients $B_i$
and so vary
holomorphically with the parameter $\lambda$.
This motivates the definition of the $\IC^*$ action in
\eqref{eqn: C* action defn downstairs},
and in turn this yields the definition of the map $\varphi$
as explained just after Lemma \ref{lem:redtrip}, by taking the projection to
$\SL_2(\IC)$ of the rank two part of the semisimplification 
of \eqref{eq: translated monod} when $\lambda=-\mu_1$.
(Note that in our conventions the $\mu_i$ are the eigenvalues of
$\sum B_i=-B_\infty$.)

\subsection*{Braid group actions}\  \  

Let us check 
that the $\IC^*$ action commutes with the braid group action on the
level of the matrices $u$.
(One suspects this is the case since the braid group actions are obtained by
integrating the isomonodromy equations, and we have seen in Proposition
\ref {prop: additive translation}
that the corresponding
$\IC$ action commutes with the Schlesinger flows.)  

The standard braid group action of the $n$-string braid group $\cB_n$
on $n$-tuples of pseudo-reflections may be given by
generators $\ga_1,\ldots,\ga_{n-1}$ with $\ga_i$ acting as
$$\ga_i(r_n,\ldots,r_1)=
(\ldots,r_{i+2},r_i,r_i^{-1}r_{i+1}r_i,r_{i-1},\ldots)$$ 
only affecting $r_i,r_{i+1}$ and preserving the product $r_n\cdots r_1$.
(For $n=3$ we previously labelled the generators differently:
$\be_1=\ga_2,\be_2=\ga_1$.)
Now suppose we write $r_i=1+e_i\otimes \al_i$ with $e_i\in V,\al_i\in V^*$,
where $V=\IC^n$.
Let us restrict to the case where the $r_i$ are linearly independent (in the
sense that any such $e_i$ form a basis of $V$).
Then it is easy to lift the above $\cB_n$ action to an action on the 
$2n^2$-dimensional  space
\beq
W:=\{\ (e_n,\ldots,e_1,\al_n,\ldots\al_1) \ \bigl|\ 
\ \{e_i\} \text{ a basis of $V$, } \al_i\in V^*, \al_i(e_i)\ne -1\ \}
\eeq
by letting $\ga_i$ fix all $e_j,\al_j$ except for $j=i,i+1$:
$$
\ga_i(\ldots,e_{i+1},e_i,\ldots,\al_{i+1},\al_i,\ldots)=
(\ldots,e_{i},r_i^{-1}e_{i+1},\ldots,\al_{i},\al_{i+1}\circ r_i,\ldots)$$
where $r_i:=1+e_i\otimes\al_i\in \GL(V)$.
(We think of $W$ as the multiplicative analogue of the space on which the
lifted equations \eqref{eq: lifted schles.} were defined.)
It is simple to check this action is well-defined on $W$.

Now we may project this lifted $\cB_n$-action to the space of the matrices $u$.
Recall the $n\times n$ matrix $u$ was defined by setting
$u_{ij}=\al_i(e_j)$.
By a straightforward computation we find that, 
if we set $u'=\ga_i(u)$ then $u'_{jk}=u_{jk}$ unless  
one of $j$ or $k$ equals $i$ or $i+1$, and
$$u'_{i i}=u_{i+1 i+1},\qquad u'_{i+1 i+1}=u_{i i}$$
$$
u'_{i i+1}=t^2_i u_{i+1 i},\qquad 
u'_{i+1 i}=u_{i i+1}/t^2_i$$
$$
u'_{i j}= u_{i+1 j}+u_{i+1 i}u_{i j},\qquad 
u'_{j i}=u_{j i+1}-u_{j i}u_{i i+1}/t_i^2
$$
$$
u'_{i+1 j}=u_{i j},\qquad 
u'_{j i+1}=u_{j i}
$$
for any $j\not\in\{i,i+1\}$ where $t^2_i:=1+u_{i i}$.

In turn $u$ contains precisely the same data as the matrix 
$$u_-^{-1}t^2u_+\in G^0\subset \GL_n(\IC)$$
where $u_\pm, t^2$ are determined by the equation
$u=t^2u_+-u_-$, and $G^0$ denotes the big-cell, consisting of the invertible
matrices that may be factorised as the 
product of a lower triangular and an upper triangular matrix.
Thus the $\cB_n$-action on $\{u\}$ is equivalent to a $\cB_n$-action on $G^0$.
Let us describe this.
First let $P_i\in \GL_n(\IC)$ denote the permutation matrix corresponding to
the permutation swapping  $i$ and $i+1$.
Thus $P_i$ equals the identity matrix except in the 
$2\times 2$ block in the $i,i+1$ position on the 
diagonal, where it equals
$\left(\begin{smallmatrix}
0&1\\ 1&0 
\end{smallmatrix}\right)$.
%
Also for any unipotent upper triangular matrix $u_+$, let 
$\xi_i(u_+)$ denote the matrix which 
equals the identity matrix except in the 
$2\times 2$ block in the $i,i+1$ position on the 
diagonal, where it equals that of $u_+$, namely
$$\left(\begin{matrix}
1&(u_+)_{i i+1}\\ 0&1 
\end{matrix}\right).$$
(This map $\xi_i$ defines a homomorphism from $U_+$ to the root group of
$\GL_n(\IC)$ 
corresponding to the $i$th simple root---cf. e.g. (3.10) \cite{bafi}.)

\begin{prop}\label{prop bn action on g0}
The induced $\cB_n$-action on $G^0$ is given by the formula
$$\ga_i(a) = P_i\xi_i(u_+)a\xi_i(u_+)^{-1}P_i.$$
where $a=u_-^{-1}t^2u_+\in G^0$. 
\end{prop}
\pf
Lifting back up to $W$, write $r_j=1+e_j\otimes\al_j$ and denote
$$(\ldots,e'_{i+1},e'_i,\ldots,\al'_{i+1},\al'_i,\ldots)=
\ga_i(\ldots,e_{i+1},e_i,\ldots,\al_{i+1},\al_i,\ldots).$$
Note that the product 
$R:=r_n\cdots r_1\in \GL(V)$ is fixed by the $\cB_n$ action.
By Theorem \ref{thm: killing} the matrix for $R$ in the 
$e_j$ basis of $V$ is $a=u_-^{-1}t^2u_+$, and in the $e'_j$ basis 
the matrix for $R$ is $\ga_i(a)$.
Thus $\ga_i(a)=S^{-1}aS$
where $S$ is the matrix for the  change of basis from $\{e_j\}$ to $\{e'_j\}$.
From the formula for the action on the $e_j$, 
$S$ equals the identity matrix except in the 
$2\times 2$ block in the $i,i+1$ position on the 
diagonal, where it equals
$$
\left(\begin{matrix}
-u_{i i+1}/t_i^2 &1 \\ 1 & 0 
\end{matrix}\right)
=
\left(\begin{matrix}
1 & u_{i i+1}/t_i^2 \\ 
0 & 1 
\end{matrix}\right)^{-1}
\left(\begin{matrix}
0& 1 \\ 
1 & 0 
\end{matrix}\right)
.$$
Finally from the equation $u=t^2u_+-u_-$ we see
$u_{i i+1}/t_i^2=(u_+)_{i i+1}$
so 
$S=\xi_i(u_+)^{-1}P_i$.
\epf

\begin{rmk}\label{rmk: qwgp}
This $\cB_n$ action on the big-cell
also  appears as the classical limit of the so-called 
quantum Weyl group actions (cf. \cite{DKP} and \cite{bafi} Remark 3.8), 
provided
we use the permutation matrices rather than Tits' extended Weyl group.
Thus we have shown, for $\GL_n(\IC)$, how the classical action of 
the quantum Weyl group is related to the standard action of $\cB_n$ on 
$n$-tuples of pseudo-reflections.
Presumably this is related to Toledano Laredo's proof \cite{VTL-duke},
for $\GL_n(\IC)$ of the
Kohno--Drinfeld theorem for quantum Weyl groups.
\end{rmk} 

\begin{cor}
The $\IC^*$ action on $\{u\}$ commutes with the $\cB_n$-action defined above.
\end{cor}
\pf
On passing to $G^0$, we recall that the $\IC^*$ action just scales $t^2$ and
leaves both $u_\pm$ fixed.
However $t^2$ does not appear in the formula of Proposition 
\ref{prop bn action on g0} for the $\cB_n$ action. 
\epf

One can now see directly why the map $\varphi$ will be $\cB_3$ equivariant.
Upon using the $\IC^*$ action to make $1$ an eigenvalue of $r_3r_2r_1$ 
we know that the $r_i$ are all block triangular in some basis.
Then we just note the obvious fact that the braid group action 
\eqref{eq: b3 on prs} on the pseudo-reflections restricts
to the action \eqref{eq: sl2 b3 action} in the $2\times 2$ block on the
diagonal. 
\end{section}

\begin{section}{Jimbo's leading term formula}\label{sn: jimbo}

So far we have described how to find some 
$\SL_2(\IC)$ triples living in
finite orbits of the braid group, and read off some properties of the
corresponding solution to Painlev\'e VI (in particular we saw
that the set of branches of
the solution correspond to the orbit under the pure braid group of the 
conjugacy classes of such triples).
In this section and the next we will describe a method to
find the corresponding solution explicitly.
This method is quite general and should work with any sufficiently generic
$\SL_2$ triple in a finite braid group orbit---in particular it is not a priori
restricted to any one-parameter family of Painlev\'e VI equations.
(One just needs to check conditions b),c),d) below for each branch of the
solution.)

The crucial step is the following formula:

\begin{thm}\ {\em(M. Jimbo \cite{Jimbo82})}
Suppose we have four matrices $M_j\in \SL_2(\IC),$ $j=0,t,1,\infty$ satisfying

a) $M_\infty M_1M_tM_0=1$,

b) $M_j$ has eigenvalues $\{\exp(\pm\pi i \theta_j)\}$ with $\theta_j\notin
\IZ$,

c) $\tr(M_0M_t)=2\cos(\pi\si)$ for some nonzero $\si\in\IC$ 
with $0\le\re(\sigma)< 1$,

d) None of the eight numbers  
$$
\theta_0\pm\theta_t\pm\sigma,\quad     \theta_0\pm\theta_t\mp\sigma, \quad
\theta_\infty\pm\theta_1\pm\sigma,\quad \theta_\infty\pm\theta_1\mp\sigma$$
is an even integer.

Then the leading term in the asymptotic expansion at zero
of the corresponding Painlev\'e VI
solution $y(t)$ 
on the branch corresponding to $[(M_0,M_t,M_1)]$ is

\beq\label{eq: JLCF}
\frac{
(\theta_0+\theta_t+\si)(-\theta_0+\theta_t+\si)(\theta_\infty+\theta_1+\si)}
{
4\si^2(\theta_\infty+\theta_1-\si)\wh s}
t^{1-\si}
\eeq
where
$$
\wh s={c\times s,}\qquad s=\frac{a+b}{d}
$$

$$a=
{e^{\pi i { \si}}}(
i\sin \left( \pi { \si} \right) \cos \left( \pi { \si_{1t}} \right) -
\cos \left( \pi { \theta_t} \right) \cos \left( \pi { \theta_\infty}\right) -
\cos \left( \pi { \theta_0}\right) \cos \left( \pi { \theta_1}\right) )
$$

$$b=i\sin \left( \pi { \si} \right) \cos \left( \pi { \si_{01}} \right)
+\cos \left( \pi { \theta_t} \right) \cos \left( \pi { \theta_1}\right) +
\cos \left( \pi { \theta_\infty} \right) \cos \left( \pi { \theta_0}\right) 
$$

{\small
$$
d=
4
\sin \left(\frac{\pi}{2}\left({ \theta_0}  +{ \theta_t}-{\si}\right)\right) 
\sin \left(\frac{\pi}{2}\left({ \theta_0}  -{ \theta_t}+{\si}\right)\right) 
\sin \left(\frac{\pi}{2}\left(\theta_\infty+{ \theta_1}-{\si}\right)\right) 
\sin \left(\frac{\pi}{2}\left(\theta_\infty-{ \theta_1}+{\si}\right)\right) 
$$
}
$$
c=\frac{
\left( \Gamma  \left( 1-{ \si} \right)  \right) ^{2}
\wh\Gamma  \left( { \theta_0}+{ \theta_t}+{ \si} \right) 
\wh\Gamma  \left(-{ \theta_0}+{ \theta_t}+{ \si} \right) 
\wh\Gamma  \left( { \theta_\infty}+{ \theta_1}+{ \si} \right) 
\wh\Gamma  \left(-{ \theta_\infty}+{ \theta_1}+{ \si}\right) 
}{
\left( \Gamma  \left( 1+{ \si} \right)  \right) ^{2}
\wh\Gamma  \left( { \theta_0}+{ \theta_t}-{ \si} \right) 
\wh\Gamma  \left(-{ \theta_0}+{ \theta_t}-{ \si} \right) 
\wh\Gamma  \left( { \theta_\infty}+{ \theta_1}-{ \si} \right) 
\wh\Gamma  \left(-{ \theta_\infty}+{ \theta_1}-{ \si}\right) 
}$$
where $\wh\Gamma(x):=\Gamma(\frac{1}{2}x+1)$ (with $\Gamma$ being the usual
gamma function) and where $\sigma_{jk}\in\IC$ (for $j,k\in\{0,t,1\}$)
is determined by
 $\tr(M_jM_k)=2\cos(\pi\sigma_{jk}),0\le\re(\sigma_{jk})\le 1$, 
so $\si=\sigma_{0t}$.
\end{thm}

\begin{rmk}
The formula \eqref{eq: JLCF} is computed directly from 
the formula \cite{Jimbo82} (2.15) for the asymptotics as $t\to 0$ 
for the coefficients
of the isomonodromic family of rank two systems.
The formulae for $\wh s$ and $s$ are as in \cite{Jimbo82} except for a sign
difference in $s$. Since this formula is crucial for us and since 
$s$ is not derived in \cite{Jimbo82} we will give a
derivation in the appendix.
\end{rmk}

\begin{rmk}
D. Guzzetti repeated Jimbo's computations in 
\cite{guzz-ellr} Section 8.3 and Appendix, but did  not reduce the formula  to
as short a form;
see \cite{guzz-ellr} (A.6) and (A.30) (but note (A.30) is not quite correct
  but is easily corrected by examining (A.28) and (A.29)). 
However we can state that the corrected version of Jimbo's formula agrees 
numerically with the corrected version of Guzzetti's
(at least for the values of the parameters used in this paper, and for those of
several hundred randomly chosen $\SL_2(\IC)$ triples). 
It is puzzling  that Guzzetti does not state that his formula does not
agree with Jimbo's.\footnote{
Also there is some confusion as to the range of validity of Jimbo's work:
namely the restriction $0\le\re(\sigma)< 1$, 
is equivalent to $\tr(M_0M_t)\not\in \IR_{\le-2}$ rather than the much
stronger condition $\vert\tr(M_0M_t)\vert\le 2$ and $\re\tr(M_0M_t)\ne-2$
appearing in \cite{guzz-ellr} (1.30).}

\end{rmk}

\begin{rmk}
To agree with Jimbo's notation we are thus relabelling the triples 
$(M_1,M_2,M_3)$ as  $(M_0,M_t,M_1)$ as well as the $\theta$ parameters
$(\theta_1,\theta_2,\theta_3,\theta_4)\mapsto
(\theta_0,\theta_t,\theta_1,\theta_\infty)$.
To keep track of this it is perhaps simplest to bear in mind the corresponding
monodromy relations
$$M_4M_3M_2M_1=1\qquad\text{and}\qquad M_\infty M_1M_tM_0=1.$$

\end{rmk}

\end{section}

\begin{section}{The Klein solution}\label{sn: klein soln}

For the Klein solution $\si$ is either $1/2$ or $1/3$ depending on the
branch.
If the solution is to be algebraic then Jimbo's formula will give 
the leading term
in the {\em Puiseux} expansion at $0$ of each branch of the solution.
Thus we find the leading term on the $j$th branch is 
of the form $C_jt^{1-\si_j}$ where
$$C_j=\frac{57}{28\wh s_j}$$
on the four branches with $\si=1/2$ ($j=0,1,2,3$) and
$$C_j=\frac{475}{308\wh s_j}$$ on the other three branches,  
having $\si=1/3$ ($j=4,5,6$).
Now we would like to evaluate these precisely on each branch and identify
them as algebraic numbers.  
A simple numerical 
inspection shows that $C_0$ has argument $\pi/4$, $C_6$ is real and
negative, 
$$C_1=-iC_0, \qquad C_2=iC_0,\qquad C_3=-C_0$$
and
$$C_4=\exp(-2\pi i/3) C_6,\qquad C_5=\exp(2\pi i/3) C_6.$$
Thus we would hope that $C_0^4$ and $C_6^3$ are rational numbers.
Using Maple we calculate the various $\wh s_j$'s numerically and then deduce:
$$C_0^4=-7/3^4, \qquad\text{so that}\qquad C_0=\frac{(1+i)7^{1/4}}{3\sqrt{2}}$$
$$C_6^3=-5^3/14, \qquad\text{so that}\qquad C_6=\frac{-5}{14^{1/3}}.$$

Thus we now know precisely the
leading coefficient $C_j$ of the Puiseux expansion at $0$ 
of each branch of the solution $y(t)$.
By substituting back into the Painlev\'e VI equation these leading terms
determine, algebraically, any desired term in the Puiseux expansion.
If the solution is to be algebraic it should satisfy an
equation of the form 
$$F(t,y(t))=0$$
for some polynomial $F(t,y)$ in two variables. 
Since the solution has $7$ branches $F$ should have degree $7$ in $y$.
Let us write $F$ in the form
$$F=q(t)y^7+p_6(t)y^6+\cdots+p_1(t)y+p_0(t)$$
for polynomials $p_i,q$ in $t$ and
define rational functions $r_i(t):=p_i/q$ for $i=0,\ldots,6$.
If $y_0,\ldots,y_6$ denote the (locally defined) solutions on the branches
then for each $t$ we have that $y_0(t),\ldots,y_6(t)$ are the roots of 
$F(t,y)=0$ and it follows that
$$y^7+r_6(t)y^6+\cdots+r_1(t)y+r_0(t)=(y-y_0(t))(y-y_1(t))\cdots(y-y_6(t)).$$
Thus, expanding the product on the right, the rational functions $r_i$ are
obtained as symmetric polynomials in the $y_i$:
$$r_0=-y_0\cdots y_6,\ \ldots\ ,r_6=-y_0-\cdots-y_6.$$ 
Since the $r_i$ are global rational functions, the Puiseux expansions of the
$y_i$ give the Laurent expansions at $0$ of the $r_i$.
Clearly only a finite number of terms of each Laurent expansion 
are required to determine each $r_i$, and indeed it is simple to convert these
truncated Laurent expansions into global rational functions.
Clearing the denominators then yields the solution curve, as in
equation 
\eqref{soln curve} of the introduction.




One may easily check on a computer that
this curve has precisely the right
monodromy over the $t$-line (and in particular is genus zero, and has monodromy
group $A_7$).
Also one finds that it has 10 singular points; 6 double points over 
$\IC\setminus\{0,1\}$ and 4 more serious
singularities over the branch points. 
Finally since it is a genus zero curve we can look for a rational
parameterisation. 
Using the CASA package \cite{CASA}, and a simple Mobius transformation, 
we find the solution may be parameterised quite simply as:

$$
y=-{\frac { \left( 5\,{s}^{2}-8\,s+5 \right)  \left( 7\,{s}^{2}-7\,s+4\right) }
{ s \left( s-2 \right)  \left( s+1 \right) \left( 2\,s-1 \right) 
\left( 4\,{s}^{2}-7\,s+7 \right)  }},\qquad
t={\frac { \left( 7\,{s}^{2}-7\,s+4 \right) ^{2}}{{s}^{3} \left( 4\,{s}^
{2}-7\,s+7 \right) ^{2}}}.$$


Note that the polynomial $F$  defining the solution curve is quite 
canonical but there are many 
possible parameterisations.
Using the parameterisation it is easy to carry out the ultimate test
and substitute back into the Painlev\'e VI equation
(with parameters $(\al,\be,\ga,\de)=(9,-4,4,45)/98$) 
finding that we do
indeed have a solution.

\end{section}

\begin{section}{Inequivalence Theorem}\label{sn: ineq}

We know (cf. Remark \ref{rmk sym. trips} above and \cite{DubMaz00} Remark 0.2)
that the five `platonic' solutions
of \cite{DubMaz00} are equivalent (via Okamoto transformations)
to solutions associated to finite subgroups of $\SL_2(\IC)$.
In other words, even though the unipotent matrices \eqref{eqn: unip} generate
an infinite group, there is an equivalent solution with finite $2\times 2$
monodromy group.

This raises the following question: 
Even though the $2\times 2$ monodromy data
\eqref{eqn: Klein 2x2}  associated to the Klein solution generates an infinite
group, is there an equivalent solution with finite $2\times 2$ monodromy?
We will prove this is not the case:

\begin{thm}\label{thm: genuine inf gp}
Suppose there is an algebraic solution of some Painlev\'e VI equation which is
equivalent to the Klein solution under Okamoto's affine $F_4$ action.
Then the corresponding $2\times 2$ monodromy data $(M_1,M_2,M_3)$ also 
generate an infinite subgroup of $\SL_2(\IC)$.
\end{thm}
\pf
First the parameters $(\theta_1,\theta_2,\theta_3,\theta_4)$ should be
equivalent to the corresponding parameters $(2,2,2,4)/7$ of the Klein
solution.
If $(M_1,M_2,M_3)$ generate a binary tetrahedral, octahedral or icosahedral
group then we will not be able to get any
sevens in the denominators (since these groups have no elements of order
seven) so any solution associated to these groups is inequivalent to the Klein
solution. (This uses the simple 
observation that Okamoto's transformations act
within the ring $\IZ[\frac{1}{2},\theta_1,\theta_2,\theta_3,\theta_4]$
so that if the $\theta_i$ are rational numbers with no sevens in the
denominators, then no equivalent set of parameters has a seven in any
denominator.) 

Next, suppose $F(t,y)=0$ is the curve defining the Klein solution and
$F_1(t_1,y_1)=0$ is the curve defining the equivalent solution.
Then we know (\cite{OkaPVI} p.361) 
that $t$, $t_1$ are related by a Mobius transformation 
permuting $0,1,\infty$.

\begin{lem}
There is an isomorphism of the curves $F(t,y)=0$ 
and $F_1(y_1,t_1)=0$ covering the automorphism of the projective line 
mapping $t$ to $t_1$.
\end{lem}
\pf
Let us recall some facts about Okamoto's transformations
from \cite{OkaPVI, NY-so8, Masuda-PVI}.
First write $q:=y, q_1:=y_1$. Then, from the formulae for the action of the
Okamoto transformations \cite{NY-so8} Table 1, \cite{Masuda-PVI} (7.14), 
we see that $q_1$ is a rational function 
of $q,p,t$, where $p$ is the conjugate variable to $q$ in the Hamiltonian
formulation of Painlev\'e VI (cf. \cite{OkaPVI} (0.6)). 
The first of
Hamilton's equations says 
$\frac{dq}{dt}=\frac{\partial H}{\partial p}$, where $H=H_{VI}$ is the
Hamiltonian \cite{OkaPVI} p.348.
By observing that $H$ is a quadratic polynomial in $p$ (and rational in $t$
and polynomial in $q$) we deduce immediately that 
$p$ is a rational function of $\frac{dq}{dt},q$ and $t$.
Moreover since $q=y$ satisfies the polynomial equation  $F(t,y)=0$ 
implicit differentiation enables us to express 
$\frac{dq}{dt}$ as a rational function of $q,t$. Thus $p$ is a rational
function of just $q, t$ and so 
in turn $q_1$ is a rational function of just $q,t$.

Now by the symmetry of the situation 
the same argument also shows $q$ is a rational function of $q_1,t_1$.
This sets up an isomorphism between the fields $\IC(q,t)\cong \IC(q_1,t_1)$
extending the isomorphism $\IC(t)\cong\IC(t_1)$ given by mapping $t$ to $t_1$.
Dualising this gives the desired isomorphism of the 
corresponding curves.
\epf

In particular we see that $F_1$ must have degree seven in $y_1$, since
the curves have the same number of branches. This implies $(M_1,M_2,M_3)$
cannot generate a cyclic group, since cyclic groups are abelian and so 
the pure
braid group acts trivially; all such solutions have just one branch.

Finally we need to rule out the binary dihedral groups which will
need more work.
Write the elements of the binary dihedral group of order $4d$ as
$$\wt{I_2(d)}=\{1,\ze,\ldots,\ze^{2d-1},\tau,\tau\ze,\ldots,\tau\ze^{2d-1}\}$$
where 
$\ze:=
\left(\begin{smallmatrix}
\varepsilon & \\
 & \varepsilon^{-1}
\end{smallmatrix}\right)
,
\tau :=
\left(\begin{smallmatrix}
0 & -1\\
1 & 0
\end{smallmatrix}\right)
$
and $\varepsilon=\exp(\pi i/d)$.

Below we will abbreviate $\tau\ze^a$ as just $\tau a$ 
and $\ze^a$ as $a$. 

The basic strategy is to go through all possible triples of elements and show
in each case that, on conjugacy classes, the generators 
$p_1:=\be^2_1,p_2:=\be^2_2$ of
the pure braid group action
cannot have two two-cycles and a three-cycle.

The basic relations we will use repeatedly  are:
$$\tau\ze^k=\ze^{-k}\tau, \qquad \ze^{2d}=1.$$
First let us record the formulae for the action of $p_1$ on all possible pairs
of elements. (Here $\be_1$ acts by mapping a pair $(x,y)$ of elements 
to $(y,y^{-1}xy)$, and $p_1$ is the square of $\be_1$.)
\begin{lem} \label{lem: pairs}
Suppose $a,b$ are arbitrary integers. Then, in abbreviated form:
$$p_1(a,b)=(a,b),$$
$$p_1(\tau a,b)=(\tau(a+2b),-b),$$ 
$$p_1(a,\tau b)=(-a,\tau(b-2a)),$$
$$p_1(\tau a, \tau b)=(\tau(2b-a),\tau(3b-2a)).$$
\end{lem}
\pf Straightforward computation.
\epf

Now, on triples, $\beta_1$ (respectively $\be_2$) maps 
$(x,y,z)$ to $(y,y^{-1}xy,z)$ (resp. $(x,z,z^{-1}yz)$).
Immediately we see that any triple of the form 
$$(a,b,c),\qquad(\tau a,b,c),\quad\text{or}\quad(a,b,\tau c)$$
will be fixed by one or both of $p_1,p_2$.
Thus the corresponding permutation
representation will have a one-cycle, which is not permitted.

In general the triples of elements of $\wt{I_2(d)}$ fall into eight `types' 
depending on if each element contains a $\tau$ or not. 
From Lemma \ref{lem: pairs} $p_1$ and $p_2$ clearly 
take triples to triples of the same type. 
After the three types already dealt with the next four are:
$$    (a,\tau b,c),\qquad(\tau a,\tau b,c),
\qquad(a,\tau b,\tau c),\qquad(\tau a, b,\tau c).$$
For each of these one finds, from Lemma \ref{lem: pairs},
that either $p_1^2$ or $p_2^2$ (or both) act trivially.
This implies that there will be no three-cycles in the permutation
representation of one or both of $p_1$ or $p_2$ 
on conjugacy classes of such triples.

Finally we need to rule out the triples of type $(\tau a,\tau b,\tau c)$.
First let us note that the conjugacy class of $\tau$ has size $d$ and
contains the elements $\tau(2a)$ for any integer $a$, 
and the conjugacy class of $\tau\ze$ has size $d$ and
contains the elements $\tau(2a+1)$.
It follows that, upto overall conjugacy, we have:
\beq \label{eq: 8th p1}
p_1(\tau a,\tau b,\tau c) \cong (\tau a,\tau b,\tau (c-k))\quad
\text{where $k:=2(b-a)$, and}
\eeq
$$p_2(\tau a,\tau b,\tau c) \cong (\tau (a-l),\tau b,\tau c)\quad
\text{where $l:=2(c-b)$.}$$
Moreover the only (possibly distinct) triple of the form 
$(\tau p,\tau b,\tau q)$ that is conjugate to $(\tau a,\tau b,\tau c)$
is $(\tau (2b-a),\tau b,\tau (2b-c))$, which is obtained by conjugating by
$\tau b$.

\begin{lem}
Let $o(k)$ be the order of the element $\ze^k$ where $k=2(b-a)$.
Then in the permutation representaion of $p_1$ the
conjugacy classes of the triple through 
$(\tau a,\tau b,\tau c)$ lies in a cycle of length $o(k)$.
\end{lem}
\pf
First if $\tau a=\tau(2b-a)$ then $\ze^k=1$ so $o(k)=1$ and 
\eqref{eq: 8th p1} says the conjugacy class of 
$(\tau a,\tau b,\tau c)$ is fixed by $p_1$. 

Secondly if $\tau a\ne\tau(2b-a)$, i.e. $o(k)>1$ then by 
\eqref{eq: 8th p1} we see
$p^r_1(\tau a,\tau b,\tau c) \cong (\tau a,\tau b,\tau (c-rk))$.
This is conjugate to $(\tau a,\tau b,\tau c)$ if and only if 
$\ze^{-rk}=1$ (using the fact that $\tau(2b-a)\ne\tau a$) 
i.e. if and only if $r$ is divisible by $o(k)$.
Thus we are in a cycle of length $o(k)$.
\epf

Similarly if $o(l)$ is the order of $\ze^l$ where $l=2(c-b)$ 
then the conjugacy class of $(\tau a,\tau b,\tau c)$ is in a cycle of
$p_2$ of length $o(l)$.

Thus in order to be equivalent to the Klein solution we need 
$o(k),o(l)\in\{2,3\}$ for all the possible $k$'s and $l$'s that occur in the
orbit. It is straightforward to check this is not possible:
First from \eqref{eq: 8th p1} note that $p_1$ maps the pair of integers 
$[k,l]$ to $[k,l-2k]$ and similarly $p_2[k,l]=[k+2l,l]$.
Thus:

i) If $o(k)=o(l)=2$ then $o(l-2k)=o(k+2l)=2$ and, repeating, we see
only two-cycles appear in the
orbit, whereas we need a three-cycle.

ii) If $o(k)=2,o(l)=3$ then $o(k+2l)=6$ and so we get an unwanted
six-cycle. (Similarly if $o(k)=3,o(l)=2$.)

iii) If $o(k)=o(l)=3$ then $\ze^{k+2l}, \ze^{l-2k}$ each 
have order either one or three. Thus 
either an unwanted one-cycle appears or we only get
three cycles; no two-cycles appear.

Thus we conclude that the Klein solution is not equivalent to any solution
coming from a finite subgroup of $\SL_2(\IC)$.
\epf
  
\end{section}

\begin{section}{Reconstruction}\label{sn: reconstruction}

Given a triple $r_1,r_2,r_3$ of generators of a three-dimensional complex
reflection group, we have explained how to obtain an $\SL_2(\IC)$ 
triple $M_1,M_2,M_3$ (in an isomorphic braid group orbit) and then how,
if Jimbo's formula is applicable, to obtain an algebraic solution $y(t)$ to
the sixth Painlev\'e equation.

In this section we will explain how to obtain from $y(t)$ a rank three
Fuchsian system with four poles on $\IP^1$ and monodromy conjugate to the
original complex reflection group (generated by three reflections).

First we recall 
(from \cite {JM81}) that the solution $y(t)$ and its derivative determine
algebraically an $\lsl_2$ system 
\begin{equation} \label{eqn: 2x2 linear recon}
\frac{d\Phi}{dz}= A(z)\Phi; \qquad A(z)=\sum_{i=1}^3\frac{A_i}{z-a_i}
\end{equation}
with monodromy $(M_1,M_2,M_3)$, where $(a_1,a_2,a_3)=(0,t,1)$, with respect to
some choice of loops generating the fundamental group of the four-punctured
sphere.  (The exact formulae will be given below.) 
Now define 
$$\wh A_i=A_i+\theta_i/2$$
for $i=1,2,3$, so that $\wh A_i$ 
has rank one (and eigenvalues $\{0,\theta_i\}$).
Then the system 
\beq\label{eq: rank two shifted conn}
\frac{d}{dz}-\sum_{i=1}^3\frac{\wh A_i}{z-a_i}
\eeq
has monodromy $(\wh M_1,\wh M_2,\wh M_3)$ where
$\wh M_i=M_i\exp(\pi\sqrt{-1}\theta_i)$, which are pseudo-reflections in
$\GL_2(\IC)$.
Write these rank one matrices as
$$\wh A_i=h_i\otimes \ga_i\quad
\text{for some $h_i\in \IC^2, \ga_i\in(\IC^2)^*,\  i=1,2,3.$}$$ 
In general the span of the $h_i$ will be two dimensional and without loss of
generality  we will suppose
that $h_1,h_2$ are linearly independent (otherwise we can relabel below). 
Now consider the three $3\times 3$ rank one matrices given by
$$
B_i:=
\left(\begin{matrix}
0 & c_i\ga_i \\
\begin{matrix} 0\\0 \end{matrix} & \wh A_i \\
\end{matrix}\right)\qquad i=1,2,3
$$
for some constants $c_1,c_2,c_3\in\IC$
and the corresponding rank three system
\beq\label{eq rk 3 syst recon}
\frac{d}{dz}-\sum_{i=1}^3\frac{B_i}{z-a_i}.
\eeq
By overall conjugation, since $h_1,h_2$ are linearly independent,
we can always assume $c_1=c_2=0$.
Now if $c_3=0$ then \eqref{eq rk 3 syst recon} is block diagonal
and reduces to \eqref{eq: rank two shifted conn}.
However if $c_3\ne 0$ we obtain a rank three system with 
$$B_i=f_i\otimes \be_i$$
for a basis $f_i$ of $V:=\IC^3$---namely:
$$
f_1=\left(\begin{matrix} 0\\h_1 \end{matrix}\right),\ 
f_2=\left(\begin{matrix} 0\\h_2 \end{matrix}\right),\ 
f_3=\left(\begin{matrix} c_3\\h_3 \end{matrix}\right),\quad
\be_i=\left(\begin{matrix} 0 & \ga_i  \end{matrix}\right).$$
Moreover, up to overall conjugation this system is independent of the choice of
nonzero $c_3$ (since conjugating by $\diag(c,1,1)$ scales $c_3$ arbitrarily).

In particular the invariant functions of the monodromy of the system 
\eqref{eq rk 3 syst recon} are independent of $c_3$ and are equal to the
invariants of the monodromy of the limiting system with $c_3=0$, since the
invariants are holomorphic functions of any parameters.

Now we can perform the scalar shift of section \ref{sn: imds}
in reverse.
Namely, in the $f_i$ basis $(B_1,B_2,B_3)$ have the form 
\beq\label{eq: Bi's in good basis}
B_1=\left(\begin{matrix}
\wt\lambda_1& b_{12} & b_{13} \\
0&0&0\\
0&0&0
\end{matrix}\right),
\qquad
B_2=\left(\begin{matrix}
0&0&0 \\
b_{21}&\wt\lambda_3&b_{23} \\
0&0&0
\end{matrix}\right),\qquad
B_3=\left(\begin{matrix}
0&0&0 \\
0&0&0 \\
b_{31}&b_{32}&\wt\lambda_3
\end{matrix}\right)
\eeq
for some numbers $b_{ij},\wt \lambda_i$.
Then the scalar shift just translates each $\wt \lambda_i$ by the same scalar.

If (as we are assuming) 
we started with a solution $y(t)$ as constructed with the procedure of
this paper then $\wt \lambda_i=\lambda_i-\mu_1$,
where $\lambda_i,\mu_i$ are related as in
\eqref{eqn: bLambda} to the original 
complex reflections $r_1,r_2,r_3$.

\begin{thm}
The system obtained 
by replacing each $\wt \lambda_i$ by $\lambda_i$ in
\eqref{eq: Bi's in good basis}
has monodromy conjugate to $(r_1,r_2,r_3)$.
In other words there is a choice of fundamental solution $\Phi$ and 
of simple positive loops
$l_i$ around $a_i$ for $i=1,2,3$ generating 
$\pi_1(\IP^1\setminus\{a_1,a_2,a_3,\infty\})$
such that $\Phi$ has monodromy $r_i$ around $l_i$.
\end{thm}

\pf
Consider the system obtained by replacing $\wt\lambda_i$ 
by $\lambda_i+\lambda$ for each $i$, for varying $\lambda$
(so $\lambda=-\mu_1$ is the original system).
Write
$$
{\bf \wh t}(\lambda)=
(t_i^2(\lambda),t_{ij}(\lambda),t_{321}(\lambda),t'_{321}(\lambda))
$$
for the invariant functions of the monodromy of the corresponding system.
These functions vary holomorphically with $\lambda$ for any 
$\lambda\in\IC$.

Write $r'_i=1+e_i\otimes\al_i, u_{ij}=\al_i(e_j)$
for the monodromy data at $\lambda=0$.
By construction the eigenvalues of $r'_1,r'_2,r'_3$ and of the product
$r'_3r'_2r'_1$ are the same as those of the $r_i$ 
(since they are determined by the residues of the Fuchsian system).

Now the invariants ${\bf \wh t}(0)$
are easily expressed in terms of 
$u$ (cf. proof of Lemma \ref{lem: 3x3 quad action})
and we know how $u$ varies with $\lambda$ (Theorem \ref{thm: FL on monod}).
It follows that 
\beq\label{eq: variation of invts}
{\bf \wh t}(\lambda)=
(h^2t_i^2,h^2t_{ij},h^2t_{321},h^4t'_{321})
\eeq
where $(t_i^2,t_{ij},t_{321},t'_{321})={\bf \wh t}(0)$
and $h=\exp(\pi i\lambda)$.

By construction we know the invariants ${\bf \wh t}(-\mu_1)$
of the original system, namely 
they equal the invariants of the  block diagonal monodromy data
$$ 
\left(\begin{matrix} 1 & \\ & \wh M_1 \end{matrix}\right),\qquad
\left(\begin{matrix} 1 & \\ & \wh M_2 \end{matrix}\right),\qquad
\left(\begin{matrix} 1 & \\ & \wh M_3 \end{matrix}\right),$$
which is the monodromy of the limiting system with $c_3=0$.

But this was set up precisely so that ${\bf \wh t}(0)$ (obtained by inverting
\eqref{eq: variation of invts} when $\lambda=-\mu_1$)
are the invariants of the original complex reflection group generators.

Finally we remark that the conjugacy class of $(r_1,r_2,r_3)$
is uniquely determined by the value 
of the invariants ${\bf \wh t}$.
This will be clear in the example below and
follows in general from the fact that the invariants ${\bf \wh t}$ 
generate the ring of
conjugation invariant functions on triples of pseudo-reflections, and that
the triple $(r_1,r_2,r_3)$ is irreducible.
\epf

\begin{rmk}\label{rmk simpler}
Having established that the resulting system has the correct monodromy,
let us record a more direct way to go from the $\wh A_i$ to the 
$B_i$ of \eqref{eq: Bi's in good basis}.
The key point is that the pairwise, and three-fold, traces of distinct
$B_i$'s are independent of the scalar shift $\lambda$ (and that the constant
$c_3$ does not contribute).
Thus if $i\ne j$ then we find
$$b_{ij}b_{ji}=\tr(B_iB_j)=\tr(\wh A_i\wh A_j),$$
$$b_{32}b_{21}b_{13}=\tr(B_3B_2B_1)=\tr(\wh A_3\wh A_2\wh A_1).$$
In general these are sufficient to determine $(B_1,B_2,B_3)$ uniquely 
up to conjugacy.
\end{rmk}

Now we will recall (from \cite{JM81})
the formulae for the $\wh A_i$ in terms 
of $y,y'$.
Let us first go in the other direction, and then invert.
Consider the following rank one matrices
$$
\wh A_1:=\left(\begin{matrix}
z_1+\theta_1 & -u z_1 \\ (z_1+\theta_1)/u & -z_1 
\end{matrix}\right),\ 
\wh A_2:=\left(\begin{matrix}
z_2+\theta_2 & -w z_2 \\ (z_2+\theta_2)/w & -z_2 
\end{matrix}\right),\ 
\wh A_3:=\left(\begin{matrix}
z_3+\theta_3 & -v z_3 \\ (z_3+\theta_3)/v & -z_3 
\end{matrix}\right)
$$
so that $\wh A_i$ has eigenvalues $\{0,\theta_i\}$ for $i=1,2,3$.
Now if we define
$$
k_1:=(\theta_4-\theta_1-\theta_2-\theta_3)/2,\qquad
k_2:=(-\theta_4-\theta_1-\theta_2-\theta_3)/2$$
and impose the equations
\beq\label{eqs: get residue at infinity}
z_1+z_2+z_3=k_2,\quad
u z_1+v z_3+w z_2=0,\quad
(z_1+\theta_1)/u+(z_3+\theta_3)/v+(z_2+\theta_2)/w=0
\eeq
then $\wh A_1+\wh A_2+\wh A_3=-\diag(k_1,k_2)$, and 
the corresponding $\lsl_2$ matrices satisfy 
$A_1+A_2+A_3=-\diag(\theta_4,-\theta_4)/2$.

Now we wish to define two $T$-invariant functions $x,y$ on 
the set of such triples $(\wh A_1,\wh A_2,\wh A_3)$, 
where the one-dimensional torus
$T\subset \SL_2(\IC)$ acts by
diagonal conjugation. 
(The function $x$ is denoted $\wt z$ in \cite{JM81}.)
First note that the $(1,2)$ matrix entry of
$$\wh A:=\frac{\wh A_1}{z}+\frac{\wh A_2}{z-t}+\frac{\wh A_3}{z-1}$$
is of the form $\frac{p(z)}{z(z-1)(z-t)}$ for some {\em linear}
polynomial $p(z)$. Thus $\wh A_{12}$ has a unique zero on the complex plane
and we define $y$ to be the position of this zero.
Explicitly one finds:
\beq\label{eqn: y}
y=\frac{t u z_1}{(t+1) u z_1+t v z_3+w z_2}.
\eeq

Then we define
\beq\label{eqn: ztilde}
x=\frac{z_1}{y}+\frac{z_2}{y-t}+\frac{z_1}{y-1}
\eeq
which is clearly $T$-invariant. 
Note that if we set 
$z=y$ then $\wh A$ is lower triangular and its first eigenvalue 
(i.e. its top-left entry) is 
$x+\frac{\theta_1}{y}+\frac{\theta_2}{y-t}+\frac{\theta_3}{y-1}$.


Now the fact is that we can go backwards and express 
the six variables 
$\{z_1,z_2,z_3,u,v,w\}$ in terms of $x,y$.
That is, given $x,y$ we wish to solve the five equations 
\eqref{eqs: get residue at infinity}, \eqref{eqn: y}, \eqref{eqn: ztilde}
in the six unknowns $\{z_1,z_2,z_3,u,v,w\}$. 
To fix up the expected one degree of freedom we impose a sixth equation
\begin{equation}\label{eq: fix torus}
{(t+1) u z_1+t v z_3+w z_2}=1
\end{equation}
so that \eqref{eqn: y} now says $y=t u z_1$. 
(This degree of freedom corresponds to the torus action mentioned above.) 
One then finds, algebraically, 
that these six equations in the six unknowns admit the unique 
solution:
{
$$z_1=
y 
\frac{
E-k_2^2 (t+1)
}{t \theta_4},\qquad
z_2=
(y-t) 
\frac{
E+
t \theta_4 (y-1){x}
k_2^2-t k_1 k_2
}{t (t-1) \theta_4}
$$$$
z_3=
-(y-1) 
\frac{
E+\theta_4 (y-t) {x}
-k_2^2 t - k_1 k_2
}{(t-1) \theta_4},
$$}
$$
u=\frac{y}{t z_1},\qquad
v=-\frac{y-1}{(t-1) z_3},\qquad
w=\frac{y-t}{t (t-1) z_2}$$
where 
$$E=
y (y-1) (y-t) {x}^2+
\bigl(\theta_3 (y-t)+t \theta_2 (y-1)-2 k_2 (y-1) (y-t)\bigr) {x}
+k_2^2 y-k_2 (\theta_3+t \theta_2).$$

Finally we recall that if $y(t)$ solves PVI then 
the variable $x$ is directly expressible in terms of
the derivative of $y$ (cf. \cite{JM81} above C55):

$$x=\frac{1}{2}\left(\frac{t (t-1) y'}{y (y-1) (y-t)}-
\frac{\theta_1}{y}-\frac{\theta_3}{y-1}-\frac{\theta_2+1}{y-t}\right).$$

Thus given a solution $y(t)$ to PVI 
we may use these formulae to reconstruct 
the matrices $\wh A_1,\wh A_2,\wh A_3$ 
upto overall conjugation by the diagonal torus.

Although it will not be needed here we remark that one needs to do a further
quadrature in order for the $\wh A_i$ to solve Schlesinger's equations;
they need to vary appropriately 
within the torus orbit. This is done via the variable
$k$ of \cite{JM81}, which evolves according to the linear
differential equation \cite{JM81} C55.

{\bf Example.}\ \ 
For the Klein solution, suppose we set the parameter $s=5/4$ (this is chosen to
give reasonably simple numbers below).
Then $t=121/125$
and the above formulae yield (cf. Remark \ref{rmk simpler}):
$$
b_{12}b_{21}=\frac{3}{224},\qquad
b_{23}b_{32}=\frac{249}{2464},\qquad
b_{13}b_{31}=\frac{5}{176},\qquad
$$
$$b_{32}b_{21}b_{13}=\frac{21}{1408}.$$
These values determine $(B_1,B_2,B_3)$  uniquely upto 
conjugacy, and it is easy to find a representative triple:
\begin{cor}
Let 
$$
B_1=\left(\begin{matrix}
\frac{1}{2}& \frac{3}{224} & \frac{21}{1408} \\
0&0&0\\
0&0&0
\end{matrix}\right),
\qquad
B_2=\left(\begin{matrix}
0&0&0 \\
1 & \frac{1}{2} & \frac{5}{176} \\
0&0&0
\end{matrix}\right),\qquad
B_3=\left(\begin{matrix}
0&0&0 \\
0&0&0 \\
\frac{332}{49} & 1 & \frac{1}{2}
\end{matrix}\right).
$$
Then the Fuchsian system
$$
\frac{d}{dz}-\left(\frac{B_1}{z}+\frac{B_2}{z-\frac{121}{125}}
+\frac{B_3}{z-1}\right) 
$$
has monodromy equal to the Klein complex reflection group,
generated by reflections.
\end{cor}

\begin{rmk}
The author is grateful to  M. van Hoeij and  J.-A. Weil
for confirming on a computer that this system does indeed admit an
invariant of degree four.
\end{rmk}

\end{section}

\begin{section}{The $3\times 3$ Fuchsian representation of PVI}
\label{sn: 3x3 rep}

In this section we will describe the direct path to the sixth Painlev\'e
equation from the $3\times 3$ isomonodromic
deformations we have been considering.
Then we will explain how to deduce a recent theorem of 
Inaba--Iwasaki--Saito \cite{IIS} from the results of this paper.

Let $V=\IC^3$ and suppose $B_1,B_2,B_3\in \End(V)$ have rank one,
linearly independent images and
$\tr(B_i)=\lambda_i$. Suppose that $B_1+B_2+B_3$ is diagonalisable with 
eigenvalues $\mu_1,\mu_2,\mu_3$ and write
$$\nabla:=d-Bdz, \qquad B(z):=\frac{B_1}{z}+\frac{B_2}{z-t}+\frac{B_3}{z-1}.$$
In section \ref{sn: imds} we showed how isomonodromic deformations of $\nabla$
lead to $\lsl_2$ isomonodromic deformations which are 
well-known to be equivalent to PVI.
(Note we are setting the pole positions $a_i$ of section \ref{sn: imds} 
to be $(a_1,a_2,a_3)=(0,t,1)$.)
One may go directly to PVI as follows (this was stated without proof in
\cite{pecr} Remark 4). 
First conjugate $B_1,B_2,B_3$ by a single element of $\GL_3(\IC)$ such that 
$B_1+B_2+B_3=\diag(\mu_1,\mu_2,\mu_3)$.
Consider the polynomial $p(z)$ defined to be the $2,3$ matrix entry of
$$z(z-1)(z-t)B(z).$$
By construction $p(z)$ is a linear polynomial, so has a unique zero on the
complex plane. Define $y$ to be the position of this zero.

\begin{prop}\label{prop: direct PVI}
If we vary $t$ and evolve $B$ according to Schlesinger's equations 
\eqref{eq: schles} then $y(t)$ solves the PVI equation 
with parameters determined by $\{\lambda_i,\mu_j\}$ as 
in \eqref{thal}, \eqref{thla}.
\end{prop}
\pf
As in section \ref{sn: imds} we perform the scalar shift by $\lambda=-\mu_1$.
Under this shift $B\mapsto B'$ say.
Then $B'_1+B'_2+B'_3=\diag(0,\mu_2-\mu_1,\mu_3-\mu_1)$ and we deduce 
$B'_i(e_1)=0$ in this basis, for each $i$. Then the bottom-right $2\times 2$
submatrix $\wh A$ of $B'$ 
also solves the Schlesinger equations and the standard theory 
\cite{JM81} says the position of the 
zero of the top-right entry of $\wh A$ solves a PVI equation.
The parameters of this PVI equation are as in \eqref{thal}, \eqref{thla}.
Clearly the position of the zero of the $1,2$ entry of $\wh A$ is the position 
of the zero of the $2,3$ entry of $B'$. It remains to check this equals that
of $B$, i.e. before the scalar shift. This is not obvious, but may be seen as
follows. 

Suppose we have performed the scalar shift by arbitrary $\lambda$.
Write $c_i(\lambda):= (B'_i)_{23}$. We claim $c_i/c_j$ is independent of
$\lambda$ for all $i,j$. This easily implies the position of the zero of
$p(z)$ is independent of $\lambda$ as required.   
To deduce the claim write $B_i=f_i\otimes \be_i$ for a basis $\{f_i\}$.
Let $\{\wh f_i\}$ be the dual basis. 
Write $\Delta:=\diag(\mu_1,\mu_2,\mu_3)$. 
Since $B_1+B_2+B_3=\Delta$ we find 
$\be_i:=\wh f_i\circ B_i=\wh f_i\circ \Delta$ 
so that $B_i= f_i\otimes \wh f_i \Delta$.
Thus upon shifting we find $B'_i= f_i\otimes \wh f_i (\Delta+\lambda)$.
Hence $c_i=\wh e_2B'_ie_3=(f_i\otimes \wh f_i)_{23}\times(\mu_3+\lambda)$ and
$c_i/c_j=(f_i\otimes \wh f_i)_{23}/(f_j\otimes \wh f_j)_{23}$
which is manifestly independent of $\lambda$. 
\epf

Of course the $2,3$ matrix entry of $B$ is not particularly special. One can
conjugate $B$ by a permutation matrix to move any of the other off-diagonal
entries of $B$ into that position.
We will show these solve equivalent Painlev\'e VI equations, related by 
Okamoto transformations. 
(Note that such conjugation corresponds to permuting the $\mu_i$.)

Let $S_3$ denote the symmetric group on three letters, 
which
is generated by the transpositions $\al=(12),\be=(13)$.
On one hand $S_3$ acts on $B$ via permutation matrices, permuting
$\mu_1,\mu_2,\mu_3$ arbitrarily.
On the other hand we may map  $S_3$ into the group of 
Okamoto transformations 
by sending
$$\al\mapsto s_1s_2s_1,\qquad\beta\mapsto (s_0s_1s_3s_4)s_2(s_0s_1s_3s_4).$$
Here each $s_i$, from \cite{NY-so8}, 
is a generator of Okamoto's affine $D_4$ symmetry group of PVI, which acts on
the set of PVI equations taking solutions to solutions.

\begin{lem}
Suppose $\si\in S_3$ and let $P$ be the corresponding permutation matrix and
let $s$ be the corresponding element of affine $D_4$.
If we vary $t$ and evolve $B$ according to Schlesinger's equations 
\eqref{eq: schles} then the PVI solution $y(t)$ associated to $PBP^{-1}$ 
is the transform, via the transformation $s$, 
of that associated to $B$.
\end{lem}
\pf
Let $y_0$ be the original solution from Proposition \ref{prop: direct PVI}.
We wish to show $y=s(y_0)$.
From $B$ and $PBP^{-1}$ we get two $2\times 2$ systems $\wh A,\wh A'$
respectively. 
It is straightforward to check their parameters are related by $s$ (cf.
Lemma \ref{lem: Okamoto action on params}).
Thus $y,s(y_0)$ solve the same PVI equation. Now
from Remark \ref{rmk simpler} we deduce that for any distinct $i,j,k$
\beq\label{eq: additive IIS}
\tr(\wh A_i\wh A_j)=\tr(\wh A'_i\wh A'_j)\qquad \text{and}\qquad
\tr(\wh A_i\wh A_j\wh A_k)=\tr(\wh A'_i\wh A'_j\wh A'_k),
\eeq
since the $P$'s cancel in the traces.
Now we point out the algebraic fact that if two solutions are related by $s$ 
then \eqref{eq: additive IIS} holds for the corresponding $2\times 2$ systems
(and generically the converse is true: \eqref{eq: additive IIS} and that the
parameters match implies the solutions match).
\epf

\begin{cor}[Inaba--Iwasaki--Saito \cite{IIS}]\label{cor: IIS thm}
Let $s$ be any element of Okamoto's affine $D_4$ symmetry group (not
necessarily in
the $S_3$ subgroup considered above).
By definition $s$ acts birationally on the space of $2\times 2$ systems $A$.
This action preserves the quadratic functions $m_{12},m_{23},m_{13}$
of the monodromy matrices of these systems.
\end{cor}
\pf
The affine $D_4$ action has five generators $s_i$ for $i=0,1,2,3,4$.
(We use the notation of \cite{NY-so8}; the indices are permuted in \cite{IIS}.)
The result is simple for $s_0,s_1,s_3,s_4$, since each of 
these transformations arise as a type of 
gauge transformation and fixes all the monodromy data
(cf. \cite{IIS} Section 6---that these cases are trivial
is stated \cite{IIS} p.15).

The hard part is to establish the result for $s_2$.
However we have shown that $s_1s_2s_1$ is the transformation which comes from
swapping $\mu_1$ and $\mu_2$.  On the level of $3\times 3$ monodromy data this
corresponds to just swapping $n_1,n_2$ (recall $n_j=\exp(\pi i \mu_j)$).
However glancing at the map $\varphi$ of \eqref{eqn: varphi} we recall
that $m_{ij}=t_{ij}/t_it_j$, which does not involve either $n_1,n_2$
and so is fixed.
\epf

\begin{rmk}
One reason we are interested 
in this result here is to check that up to equivalence
there is just one Klein solution. Recall \cite{pecr} that there are exactly
{\em two} braid group orbits of conjugacy classes of triples of generating
reflections of the Klein complex reflection group, containing the triples
$$(r_1,r_2,r_3)\qquad\text{and}\qquad (r_3,r_2,r_1)$$
respectively. It is easy to see they have PVI parameters which are equivalent
under the affine $D_4$ group
(and have isomorphic $\cP_3$ orbits).
But one would like to check the actual PVI solutions are equivalent.
This is facilitated by Corollary \ref{cor: IIS thm}: it is sufficient 
to check the quadratic functions of the corresponding $2\times 2$ 
monodromy data match up.  
In turn, via $\varphi$, this amounts to checking that
$\tr(r_1r_2)=\tr(r_3r_2)$. But we saw in section \ref{sn: braid orbits}
both sides equal $1$. 
\end{rmk}
\end{section}

\ 

\appendix
\begin{section}{}\label{apx: bjl extraction}

We wish to explain how to extract Theorem 
\ref{thm: FL on monod} from the paper \cite{BJL81} of Balser, Jurkat and Lutz.

Let us briefly recall the set-up of \cite{BJL81}.
Given an $n\times n$ matrix $A_1$ and a diagonal matrix
$B_0=\diag(b_1,\ldots,b_n)$  (with the $b_i$ pairwise distinct)
one considers the Fuchsian connection (\cite{BJL81} (0.2))
\beq \label{eqn: fuchsian bjl connn}
d-(B_0-z)^{-1}(1+A_1)dz
=d-\sum_{i=1}^n \frac{-E_i(1+A_1)}{z-b_i}dz
\eeq
where $E_i$ is the $n\times n$ matrix with a one in its $i,i$ entry and is
otherwise zero. 
(To avoid confusion with other notation of the present paper we have
relabelled $t\mapsto z, \Lambda\mapsto B_0, 
\lambda_i\mapsto b_i$ from \cite{BJL81}.)
Write $\Lambda'=\diag(\lambda_1',\ldots,\lambda_n')$ for the diagonal part of
$A_1$ and suppose that 

$i)$ No $\lambda_i'$ is an integer, and that 

$ii)$ No eigenvalue of $A_1$ is an integer.

Condition $i)$ implies that each residue of \eqref{eqn: fuchsian bjl connn}
is rank one and has non-integral trace.

Now one chooses an admissible branch cut direction $\eta$ 
and cuts the complex $z$-plane from each $b_i$ to $\infty$ along the direction 
$\eta$, leaving a simply connected domain $\cP_\eta\subset\IC$.
(In fact (\cite{BJL81} p.694) one takes $\eta\in\IR$ and uses $\eta$ to 
give logarithm choices on $\cP_\eta$ near each $b_i$.)
The direction $\eta$ is said to be `admissible' if none of these cuts overlap,
and the inadmissible $\eta$ in the interval $(-\pi/2,3\pi/2]$ are labelled
$\eta_0,\ldots,\eta_{m-1}$ (with $\eta_{i+1}<\eta_i$).
This labelling is extended to all integral subscripts $\nu$ by setting
$\eta_{\nu}=\eta_{\nu+km}+2\pi k$ for any integer $k$.

Now, given an admissible $\eta$ one may canonically construct a certain
fundamental solution $Y^*(z)$ of \eqref{eqn: fuchsian bjl connn} on $\cP_\eta$
and define an $n\times n$ matrix $C=C(\eta)$, with ones on the diagonal, 
such that (by \cite{BJL81} Lemma 1):
After continuing $Y^*$ along a small positive loop around $b_k$ (and crossing
the $k$th cut), $Y^*$ becomes 
$$Y^*(1+C_k^*)$$
where $C_k^*$ is zero except for its $k$th column which equals the $k$th
column of $C\wt D$, where $\wt D:=\exp(-2\pi i \Lambda')-1$.

If $\eta$ varies through admissible values then $C(\eta)$ does not change.
Thus we choose an integer $\nu$ and let $C=C_\nu$ be $C(\eta)$ for any 
$\eta\in(\eta_{\nu+1},\eta_\nu)$, as in \cite{BJL81} Remark 3.1 p.699.

Thus if we suppose (\cite{BJL81} p.697) that, when looking along $\eta$
towards infinity, that $b_{k+1}$ lies to the right of $b_k$ for 
$k=1,\ldots,n-1$, then the monodromy of $Y^*$ around a large positive loop
encircling all the $b_i$ is the product of pseudo-reflections
$$(1+C_1^*)\cdots(1+C_n^*).$$

The main facts we need from \cite{BJL81} now are:

1) That the Stokes matrices $C_\nu^\pm$ of the irregular connection
$$d-\left(B_0 + \frac{A_1}{x}\right)dx$$
are determined by $C_\nu$ by the equation
\beq\label{eqn: CCpCm}
C_\nu D = C^+_\nu-e^{2\pi i \Lambda'}C^-_\nu
\eeq
where 
$D:=1-\exp 2\pi i\Lambda'$.
(This is equation (3.25) of \cite{BJL81}, and that the $C_\nu^\pm$ are the Stokes
matrices is the content of \cite{BJL81} Theorem 2, p.714.)

2) That the Stokes matrices are unchanged if $A_1$ undergoes a scalar shift
$A_1\mapsto A_1-\lambda$ (\cite{BJL81} Remark 4.4, p.712).

This is sufficient to determine how the pseudo-reflections 
$1+C_k^*$ vary under the scalar shift; one doesn't need to know how the Stokes
matrices are defined, only that they are triangular matrices with ones on the
diagonal.
The main subtlety one needs to appreciate is that:
{\em in the above convention (with $b_{k+1}$ to the right of $b_k$)
then 
\beq\label{eq cpm op types}
C^+_\nu \text{ {is lower triangular and }}
C^-_\nu \text{ {is upper triangular}}.
\eeq
}
Indeed (\cite{BJL81} p.701, paragraph before (3.16)) 
$C^{+/-}_\nu$ are upper/lower triangular respectively
if $b_1,\ldots,b_n$ are ordered according to the dominance relation
on $S'_{\nu+1}$.
This dominance relation is defined (top of p.699) so that
it coincides with the natural ordering of the indices
if $b_1,\ldots,b_n$ are ordered so that the $j$th cut (along the direction
$\eta\in(\eta_{\nu+1}, \eta_\nu)$)
lies to the right of the $k$th cut whenever $j<k$ (again looking along $\eta$
towards infinity).
This is {\em opposite} to the previous ordering of the $b_i$.
Thus sticking to our original ordering we deduce \eqref{eq cpm op types}.

In summary if we set $V=\IC^n$ and write 
$1+C^*_k=1+v_k\otimes \be_k$ (where $\{\be_i\}$ is the standard basis of $V^*$
and $v_k\in V$ is the $k$th column of $v:=C_\nu \wt D$, so 
$v_{ij}=\be_i(v_j)$) then 
$$v=C_\nu\wt D = (C_\nu^+-e^{2\pi i \Lambda'}C_\nu^-)D^{-1}\wt D\sim
e^{-2\pi i \Lambda'}C_\nu^+-C_\nu^-$$
where we note that $e^{2\pi i \Lambda'}\wt D=D$ and where
$\sim$ is defined so that $A\sim B$ if there is an invertible 
diagonal matrix $s$ such that $A=sBs^{-1}$. 
(This conjugation by $s$ just corresponds to different choices of 
$v_k,\be_k$ such that $1+C^*_k=1+v_k\otimes \be_k$ and so 
clearly does not affect the corresponding pseudo-reflections.)

Thus under the scalar shift, the upper triangular part of 
$v$ is fixed, the lower triangular part is scaled by $\exp(2\pi i \lambda)$
and the diagonal part 
$e^{-2\pi i \Lambda'}-1$ is changed to $e^{-2\pi i (\Lambda'-\lambda)}-1$.

Finally let us relate this back to our conventions in the body of the paper.
Namely we have a connection
$$d-\sum\frac{B_i}{z-a_i}dz$$
with rank one residues, monodromy $r_i$ around $a_i$ and monodromy 
$r_n\cdots r_1$ around a large positive loop.
The images of the $B_i$ make up a basis of $V$ so we may conjugate 
$(B_1,\ldots, B_n)$ such that each $B_i$ is zero except in row $(n-i+1)$.

Then we set 
$b_i=a_{n-i+1}$ and define  $A_1=-1-\sum B_i$ 
so that 
$$\frac{-E_i(1+A_1)}{z-b_i}=\frac{B_{n-i+1}}{z-a_{n-i+1}}$$
for each $i$ and that
$$(1+C^*_1,\ldots,1+C^*_n)=(r_n,\ldots, r_1)$$
upto overall conjugation.
Thus if we write 
$r_i=1+e_i\otimes \al_i$ and define $u$ by $u_{ij}=\al_i(e_j)$
we have that
$$u\sim\Omega v \Omega$$
where $\Omega$ is the order reversing permutation matrix 
($\Omega_{ij}=\delta_{i n-j+1}$).
Note that we denoted the trace of $B_i$ as $\lambda_i$ 
so that the diagonal part $\Lambda'$ of $A_1$ is
$$\Lambda'=-1-\Omega\Lambda\Omega$$
where $\Lambda:=\diag(\lambda_1,\ldots,\lambda_n)$.
Thus, if we write
$$u=t^2u_+-u_-$$
with $t^2$ diagonal and $u_{+/-}$ upper/lower triangular with ones on the
diagonal, we have that, upto overall conjugation by a diagonal matrix:
$$u_+=\Omega C_\nu^+\Omega,\qquad u_-=\Omega C_\nu^-\Omega$$
and $t^2=\exp(2\pi i \Lambda)$.
Therefore under the scalar shift both $u_\pm$ are fixed and so
$u$ is changed to $h^2t^2u_+-u_-$, establishing 
Theorem \ref{thm: FL on monod}.

\begin{rmk} \label{rmk: up/low distinction}
One can check independently that it is the upper triangular part of $u$
that should be scaled by $h^2$, rather than the lower triangular part, 
since we know that the eigenvalues of $r_n\cdots r_1$ should be scaled
by $h^2$. Indeed if we expand
$$\tr(r_n\cdots r_1)=
n+\sum_i u_{ii}+\sum_{i>j}u_{ij}u_{ji}+
\sum_{i>j>k}u_{ij}u_{jk}u_{ki}+\cdots+
u_{nn-1}u_{n-1n-2}\cdots u_{21}u_{1n}$$
$$=\sum_i t_i+\sum_{i>j}u_{ij}u_{ji}+
\sum_{i>j>k}u_{ij}u_{jk}u_{ki}+\cdots+
u_{nn-1}u_{n-1n-2}\cdots u_{21}u_{1n}$$
we see that scaling just the upper triangular part of $u$ (and 
the $t_i$) scales each term here by $h^2$ as required, and
otherwise one obtains higher powers of $h$.
\end{rmk}

\end{section}

\begin{section}{}\label{apx: jimbo's s}

We 
will explain how to derive the formula for the parameter $s$
in Jimbo's formula \eqref{eq: JLCF}.
This formula is stated incorrectly, and not derived, in \cite{Jimbo82}.
Since it is not immediately clear how to derive the formula we sketch the main
steps here, and point out the (probably typographical) error.
(We remark that the whole procedure described in the present 
paper does not work without this correction.)

Suppose we have four matrices $M_j\in \SL_2(\IC),$ $j=0,t,1,\infty$ satisfying
\begin{equation} \label{eq: prod. reln}
M_\infty M_1M_tM_0=1,
\end{equation}
and $M_j$ has eigenvalues $\{\exp(\pm\pi i \theta_j)\}$ where $\theta_j\notin
\IZ$. 
Write $\varepsilon_\infty=\exp(\pi i \theta_\infty)$ and 
suppose $M_\infty$ is actually diagonal ($M_\infty=\diag(\epi,\epi^{-1})$).
Define $\sigma_{jk}\in\IC$ with $0\le\re(\sigma_{jk})\le 1$ 
(for $j,k\in\{0,t,1,\infty\}$) by
$$\tr(M_jM_k)=2\cos(\pi\sigma_{jk}),$$
and let $\sigma:=\sigma_{0t}$ and $\eps:=\exp(\pi i \si)$.

Under the further assumptions that $\si$ is nonzero, 
that $0\le\re(\sigma)< 1$ and that
none of the eight numbers
$$
\theta_0\pm\theta_t\pm\sigma,\quad     \theta_0\pm\theta_t\mp\sigma, \quad
\theta_\infty\pm\theta_1\pm\sigma,\quad \theta_\infty\pm\theta_1\mp\sigma$$
is an even integer,
Jimbo \cite{Jimbo82} p.1141 points out that, up to overall conjugacy by a
diagonal matrix, $M_0,M_t,M_1$ are given, for some $s\in\IC^*$, by:
$$
(i s_\si) M_0=
C^{-1}\left(
\begin{matrix}
\eps c_0-c_t & 2 s\al'\ga' \\
-2s^{-1}\be'\de' & -\epsi c_0+ c_t
\end{matrix}
\right)C,
\quad
(i s_\si) M_t=
C^{-1}\left(
\begin{matrix}
\eps c_t-c_0 & -2 s\eps\al'\ga' \\
2s^{-1}\epsi\be'\de' & -\epsi c_t+ c_0
\end{matrix}
\right)C
$$
$$
(i s_\infty) M_1=
\left(
\begin{matrix}
c_\sigma-\epii c_1 & -2\epii\be\ga \\
2\epi\al\de & -c_\sigma+\epi c_1
\end{matrix}
\right),
\qquad\text{where}\qquad 
C=\left(
\begin{matrix}
\de & \be\\
\al & \ga
\end{matrix}
\right)
$$
and we have used the temporary notation:
$$c_j:=\cos(\pi\theta_j),\quad c_\si:=\cos(\pi\si),\quad
 s_j:=\sin(\pi\theta_j), \quad s_\si:=\sin(\pi\si),$$
$$\al=\sin \frac{\pi}{2}(\theta_\infty-\theta_1+\sigma),
\be=\sin \frac{\pi}{2}(\theta_\infty+\theta_1+\sigma),
\ga=\sin \frac{\pi}{2}(\theta_\infty+\theta_1-\sigma),
\de=\sin \frac{\pi}{2}(\theta_\infty-\theta_1-\sigma),$$
$$
\al'=\sin \frac{\pi}{2}(\theta_0-\theta_t+\sigma),
\be'=\sin \frac{\pi}{2}(\theta_0+\theta_t+\sigma),
\ga'=\sin \frac{\pi}{2}(\theta_0+\theta_t-\sigma),
\de'=\sin \frac{\pi}{2}(\theta_0-\theta_t-\sigma).$$

Notice that $\sigma_{01}$ and $\sigma_{1t}$ do not appear in these formulae;
The idea now is to express the parameter $s$ in terms of 
 $\sigma_{01}, \sigma_{1t}$ (and the other parameters). 
Naively we can just calculate $\tr(M_0M_1)$ and $\tr(M_1M_t)$ 
from the above formulae, equate them with
$2\cos(\pi\sigma_{01}), 2\cos(\pi\sigma_{1t})$ respectively and try to solve
for $s$. However this yields a complicated expression in the trigonometric
functions and a simple formula looks beyond reach.

The key observation to simplify the computation is that the above
parameterisation of the matrices is such that
$CM_tM_0C^{-1}$ is diagonal and equal to $\Delta:=\diag(\eps,\epsi)$ (and also
equal to  
$CM_1^{-1}M_\infty^{-1}C^{-1}$ by \eqref{eq: prod. reln}).
Thus we find 
\begin{equation}\label{eq: 10 coseqn}
2\cos(\pi\sigma_{01})= \tr(M_0M_1)=
\tr\bigl(CM_\infty^{-1}C^{-1}\Delta^{-1}(CM_0C^{-1})\bigr)
\end{equation}
\begin{equation}\label{eq: 1t coseqn}
2\cos(\pi\sigma_{1t})= \tr(M_1M_t)=
\tr\bigl(CM_\infty^{-1}C^{-1}\Delta^{-1}(CM_tC^{-1})\bigr)
\end{equation}
whose right-hand sides are more manageable expressions 
in the trigonometric functions, 
and are linear in $1,s,s^{-1}$.
If we take the combination \eqref{eq: 10 coseqn}$+\eps$\eqref{eq: 1t coseqn}
then the $s^{-1}$ terms cancel and upon rearranging we find:
\begin{align} \label{eq: first s formula}
2i\det(C)s_\si(c_{01}+\eps c_{1t})=
&(\ga\de\epii-\al\be\epi)(\eps-\epsi)c_2+\notag \\
&(\ga\de\epi-\al\be\epii)(\eps^2-1)c_0+
2s\al\ga\al'\ga'(\epi-\epii)(\eps-\epsi)
\end{align}
where $c_{01}=\cos(\pi \si_{01})$ and $c_{1t}=\cos(\pi \si_{1t})$.
To proceed we note:

\begin{lem}\ 

a) $\ga\de-\al\be:=\det(C)=-s_\infty s_\si$

b) $\ga\de\epii-\al\be\epi=is_\infty(\eps c_\infty-c_1)$

c) $\ga\de\epi-\al\be\epii=is_\infty(c_1-\epsi c_\infty)$
\end{lem}
\pf
A few applications of  standard trigonometric  formulae yield
$$\al\be=(c_1-(c_\infty c_\si-s_\infty s_\si))/2,\qquad
\ga\de=(c_1-(c_\infty c_\si+s_\infty s_\si))/2$$
which gives $a)$ immediately and also enable $b), c)$  to be easily deduced.
\epf

Substituting these into \eqref{eq: first s formula} and cancelling a factor
of $2s_\si s_\infty=-(\eps-\epsi)(\epi-\epii)/2$ yields the desired formula:

$$s=\frac{
\eps(i s_\si c_{1t}-c_tc_\infty-c_0c_1)+
is_\si c_{01}+c_tc_1 +c_\infty c_0}{4\al\ga\al'\ga'}.
$$
This differs from formula $(1.8)$ of \cite{Jimbo82} in a single sign, in
$\al$, which was crucial for us. 
\end{section}

\renewcommand{\baselinestretch}{1}              %
\normalsize
\bibliographystyle{amsplain}    \label{biby}
\bibliography{../thesis/syr}    
\end{document}